\documentclass[11pt]{gtart}

\usepackage{amsmath,amssymb}
\usepackage{graphics}
\usepackage{psfrag}
\usepackage{amscd}

\newcommand{\Exp}[2]{\ensuremath{\exp_{#1}\!{#2}}}
\newcommand{\s}{\ensuremath{S^1}}
\newcommand{\exps}[1]{\Exp{#1}{\s}}

\newcommand{\expg}[2]{\Exp{#1}{\Gamma_{#2}}}
\newcommand{\expgb}[2]{\Exp{#1}{(\Gamma_{#2},v)}}
\newcommand{\expinf}[1]{\Exp{}{\, {#1}}}
\newcommand{\expsig}[2][{}]{\Exp{#2}{\Sigma_{#1}}}
\newcommand{\sym}[2]{\ensuremath{\operatorname{Sym}^{#1}\!{#2}}}
\newcommand{\symthree}{\sym{3}{\Sigma}}
\newcommand{\symtwo}{\sym{2}{\Sigma}}
\newcommand{\stwo}{\ensuremath{S^2}}
\newcommand{\expstwo}[1]{\Exp{#1}{\stwo}}
\newcommand{\expstwob}[1]{\Exp{#1}{(\stwo,\ast)}}

\newcommand{\twobytwo}[4]
{\left[\begin{array}{cc} #1 & #2 \\ #3 & #4 \end{array}\right]}
\newcommand{\sbinom}[2]{{\genfrac{[}{]}{0pt}{}{#1}{#2}}_{-1}}
\newcommand{\qbinom}[2]{\genfrac{[}{]}{0pt}{}{#1}{#2}}
\newcommand{\drop}{\partial}

\newcommand{\vcell}[2][{}]{\ensuremath{\sigma^{#2}_{#1}}}
\newcommand{\ecell}[2][{}]{\ensuremath{\tilde\sigma_{#1}^{#2}}}
\newcommand{\fcell}[2][{}]{\ensuremath{\tilde\tau_{#1}^{#2}}}
\newcommand{\vfcell}[2][{}]{\ensuremath{\tau_{#1}^{#2}}}
\newcommand{\qvfcell}[1]{\ensuremath{\tau_{\rational}^{#1}}}
\newcommand{\dlambda}{\partial^\lambda}
\newcommand{\dnu}{\partial^\nu}
\newcommand{\dgamma}{\partial^\gamma}
\newcommand{\dq}{\partial^\lambda_\rational}
\newcommand{\db}{\bar{\partial}}

\newcommand{\decarrow}[1]{\stackrel{#1}{\longrightarrow}}
\newcommand{\tH}{\tilde{H}} 
\newcommand{\one}{\mathsf{1}}
\newcommand{\two}[1]{\mathsf{2}_{#1}}
\newcommand{\ttwo}[1]{\tilde{\mathsf{2}}_{#1}}
\newcommand{\three}[1]{\mathsf{3}_{#1}}
\newcommand{\four}[1]{\mathsf{4}_{#1}}

\def\co{\colon\thinspace}

\renewcommand{\S}{\mathcal{S}}        
\newcommand{\E}{\mathcal{E}}        

\newcommand{\real}{\ensuremath{\mathbf{R}}}
\newcommand{\integer}{\ensuremath{\mathbf{Z}}}
\newcommand{\rational}{\ensuremath{\mathbf{Q}}}
\newcommand{\complex}{\ensuremath{\mathbf{C}}}

\newcommand{\branch}{\ensuremath{\mathcal{D}}}
\newcommand{\model}{\ensuremath{E_3\Sigma}}

\DeclareMathOperator{\id}{id}
\DeclareMathOperator{\coker}{coker}
\DeclareMathOperator{\intr}{int}
\DeclareMathOperator{\Span}{span}
\DeclareMathOperator{\sign}{sign}
\DeclareMathOperator{\rank}{rank}
\renewcommand{\int}{\intr}

\newtheorem{theorem}{Theorem}
\newtheorem{lemma}{Lemma}

\numberwithin{equation}{section}

\begin{document}

\title{Finite subset spaces of closed surfaces}
\author{Christopher Tuffley}
\address{Department of Mathematics, University of California at Davis \\
         One Shields Avenue, Davis, CA 95616-8633, U.S.A.}
\email{tuffley@math.ucdavis.edu}

\begin{abstract}
The $k$th finite subset space of a topological space $X$ is the space
\Exp{k}{X}\ of non-empty finite subsets of $X$ of size at most $k$,
topologised as a quotient of $X^k$. The construction is a homotopy
functor and may be regarded as a union of configuration spaces of
distinct unordered points in $X$.
We show that the finite subset spaces of a connected $2$--complex admit
``lexicographic cell structures'' based on the lexicographic order on
$I^2$ and use these to study the finite subset spaces of closed
surfaces. We completely calculate the rational homology of the finite
subset spaces of the two-sphere, and determine the top integral homology
groups of $\expsig{k}$ for each $k$ and closed surface $\Sigma$. In
addition, we use Mayer-Vietoris arguments and the ring structure
of $H^*(\sym{k}{\Sigma})$ to calculate the integer cohomology groups
of the third finite subset space of $\Sigma$ closed and orientable.
\end{abstract}


\maketitle

\section{Introduction}

The $k$th finite subset space of a topological space $X$ is the space
$\Exp{k}{X}$ of nonempty subsets of $X$ of size at most $k$,
topologised as a quotient of $X^k$ via the map sending each $k$--tuple
to the set consisting of its entries.
The construction is a homotopy functor, and if $X$ is compact \Exp{k}{X}\
may be regarded as a compactification of the configuration space
of unordered $k$--tuples of distinct points in $X$.

In our previous papers~\cite{circles02,graphs03} we studied the finite
subset spaces of the circle, connected graphs, and punctured surfaces.
In this sequel we take the first steps towards 
an understanding of the
finite subset spaces of closed surfaces. We do this in two 
orthogonal directions. We first use ideas and techniques developed 
in~\cite{circles02,graphs03} to show that finite subset 
spaces of connected
$2$--complexes admit ``lexicographic cell structures'' based
on the lexicographic ordering of $I^2$. Applying these to the
standard cell structures for closed surfaces we completely calculate
the rational homology of the finite subset spaces of $S^2$, and 
determine the top integral homology groups of $\expsig{k}$
for each $k$ and closed surface $\Sigma$. We then use a 
quite different approach to completely calculate the integral
cohomology groups of the third finite subset space of a closed
orientable surface. This is the first nontrivial case,
since the second finite subset
space co-incides with the second symmetric product, for which the answer
is already known. We build a homotopy model for \expsig{3}\ out
of \sym{3}{\Sigma}\ and the mapping cylinder of $\Sigma^2\rightarrow
\sym{2}{\Sigma}$, and calculate the cohomology using 
the Mayer-Vietoris sequence and the calculation
of $H^*(\sym{k}{\Sigma})$ due to Macdonald~\cite{macdonald62} and
Seroul~\cite{seroul72-pdm,seroul72-crasci}.

\subsection{Finite subset spaces}

We begin by recalling some basic facts and constructions on finite subset
spaces. For a history and bibliography see~\cite{graphs03}. 

The $k$th finite subset space of $X$ may be viewed as the quotient of 
the $k$th symmetric product $\sym{k}{X}=X^k/S_k$ obtained by forgetting 
multiplicities. The constructions co-incide for $k=1$ and $2$ but
differ for $k\geq 3$, since the points $(a,a,b)$ and $(a,b,b)$ in 
$X^3$ map to distinct points in \sym{3}{X}\ but the same point in 
\Exp{3}{X}. In particular $\Exp{1}{X}=X$, $\Exp{2}{X}=\sym{2}{X}$, and
\Exp{k}{X}\ is a proper quotient of $\sym{k}{X}$ for $k\geq 3$.

A second view of \Exp{k}{X}\ is obtained by regarding it as a union of 
configuration spaces of unordered tuples of distinct points in $X$.
For $j\leq k$ there is a natural inclusion map
\begin{equation}
\Exp{j}{X}\hookrightarrow\Exp{k}{X} : \Lambda \mapsto \Lambda,
\label{inclusion.eq}
\end{equation}
and if $X$ is Hausdorff this is a homeomorphism onto its image~\cite{handel00}.
In this case each stratum $\Exp{j}{X}\setminus\Exp{j-1}{X}$ is homeomorphic
to the configuration space of unordered $j$--tuples of distinct points in
$X$, so that \Exp{k}{X}\ may be regarded as a union of these spaces, with
the topology recognising that configurations of different cardinalities may
be considered close. We define the full finite subset space \expinf{X}\ 
to be the 
direct limit of the system~\eqref{inclusion.eq} of inclusions,
\[
\expinf{X} = \varinjlim\Exp{k}{X}.
\]
Both \Exp{k}{}\ and \expinf{}\ may be made into functors in the obvious
way: given $f\co X\rightarrow Y$ we define $\Exp{k}{f}$ and $\expinf{f}$
by sending $\Lambda\subseteq X$ to $f(\Lambda)\subseteq Y$. The 
homotopy classes of \Exp{k}{f}\ and \expinf{f}\
depend only on the homotopy class of $f$, making both \Exp{k}{}\ and 
\expinf{}\ homotopy functors.

Given a basepoint $x_0\in X$ we define the $k$th based finite subset
space to be the subspace 
\[
\Exp{k}{(X,x_0)} = \{\Lambda\in\Exp{k}{X}| x_0\in \Lambda\}.
\]
This subspace is the image of the map $\cup\{x_0\}$ 
that adjoins $x_0$ to each element of \Exp{k-1}{X},
\[
\cup\{x_0\}\co\Exp{k-1}{X} \rightarrow \Exp{k}{X} : 
\Lambda\mapsto \Lambda\cup\{x_0\}.
\]
It should be noted that \Exp{k-1}{X}\ and \Exp{k}{(X,x_0)}\ are in
general topologically different, as $\cup\{x_0\}$ is generically 
two-to-one on the subspace $\Exp{k-2}{X}$ of $\Exp{k-1}{X}$. The based
finite subset spaces \Exp{k}{(X,x_0)}\ are often more tractable than
the unbased spaces, and frequently play an important role as stepping
stones to understanding them.

For each $k$ and $\ell$ the isomorphism 
$X^k\times X^\ell\rightarrow X^{k+\ell}$ descends to a map
\[
\cup\co\Exp{k}{X}\times\Exp{\ell}{X}\rightarrow\Exp{k+\ell}{X}
\]
sending $(\Lambda_1,\Lambda_2)$ to $\Lambda_1\cup\Lambda_2$. This leads to a 
form of product on
maps $g\co Y\rightarrow \Exp{k}{X}$, $h\co Z\rightarrow \Exp{\ell}{X}$, and we
define $g\cup h\co Y\times Z\rightarrow \Exp{k+\ell}{X}$ to be the composition
\[
Y\times Z \xrightarrow{g\times h} \Exp{k}{X}\times\Exp{\ell}{X} 
\xrightarrow{\cup} \Exp{k+\ell}{X}.
\]
We will make extensive use of this product in constructing our cell
structures for finite subset spaces, and we note that 
$(f\cup g)\cup h = f\cup (g\cup h)$.

\subsection{Summary of main results}

In this section and elsewhere in the paper we adopt the convention
that where a co-efficient group or ring is not specified integer
co-efficients should be assumed.

We first show that the finite subset spaces of a connected finite
$2$--complex admit lexicographic cell structures in 
section~\ref{lexicographic.sec},
and then we use these in section~\ref{twosphere.sec} to completely 
calculate the rational homology of $\expstwob{k}$ and \expstwo{k}.
\begin{theorem}
The space $\expstwob{k}$ has the rational homology of $S^{2k-2}$, and
the space $\expstwo{k}$ has the rational homology of $S^{2k}\vee S^{2k-2}$.
\label{s2rational.th}
\end{theorem}

\noindent
More careful attention at the top end of the chain complex shows that
``rational'' cannot be replaced by ``integral'' in this theorem for
$k\geq 4$ in the based case and $k\geq 3$ in the unbased case:
\begin{theorem}
The top three integral homology groups of $\expstwob{k}$ are
$\integer$ in dimension $2k-2$, $\{0\}$ in dimension $2k-3$, and
$\integer/(k-2)\integer$ in dimension $2k-4$. The group
$H_{2k-4}$ is generated by the top homology class of
$\expstwob{k-1}\hookrightarrow\expstwob{k}$. 

The top three integral homology groups of \expstwo{k}\ are
$\integer$ in dimension $2k$, $\{0\}$ in dimension $2k-1$, and
$\integer\oplus\integer/(k-1)\integer$ in dimension $2k-2$.
The group $H_{2k-2}$ is generated by the top classes
$[\expstwo{k-1}]$ of \expstwo{k-1}\ and $[\expstwob{k}]$ of
\expstwob{k}, subject to the relation 
$(k-1)\bigl([\expstwo{k-1}]-2[\expstwob{k}]\bigr)=0$. 
\label{s2integral.th}
\end{theorem}

\noindent
In addition, we calculate the integral homology completely as far
as $k=6$ in the based case and $k=5$ in the unbased case.

We then turn to the finite subset spaces of higher genus surfaces in
section~\ref{highergenus.sec}. 
Handel~\cite{handel00} has shown that 
$H^{nk}(\Exp{k}{M^n};\integer/2\integer)$
has rank one for a closed connected $n$--manifold, $n\geq 2$, and using the
cell structures 
constructed in section~\ref{lexicographic.sec} we
prove the following refinement of this result for $n=2$:
\begin{theorem}
\label{surfacehomology.th}
Let $\Sigma$ be a closed surface of genus $g$. Then
\[
H_{2k}(\expsig{k})=\begin{cases}
                    \integer & \mbox{if $\Sigma$ is orientable,} \\
                    0    & \mbox{if $\Sigma$ is non-orientable,}
                   \end{cases}
\]
and
\[
H_{2k-1}(\expsig{k})=\begin{cases}
                    \integer^{2g} & \mbox{if $\Sigma$ is orientable,} \\
                    \integer/2\integer & \mbox{if $\Sigma$ is non-orientable.}
                    \end{cases}
\]
\end{theorem}

Given an orientation on $\Sigma$ we may canonically orient the
manifold $\expsig{k}\setminus\expsig{k-1}$ by orienting each
tangent space as the direct sum
\begin{equation}
T_\Lambda(\expsig{k}\setminus\expsig{k-1}) = 
        \bigoplus_{p\in\Lambda} T_p\Sigma .
\label{surfaceorientation.eq}
\end{equation}
Since $\Sigma$ is even dimensional the order of the summands does not
matter and we obtain a consistent orientation on the configuration
space $\expsig{k}\setminus\expsig{k-1}$. This in turn orients $H_{2k}$,
and Theorem~\ref{surfacehomology.th} then implies that a map between the
$k$th finite subset spaces of two closed oriented surfaces has
a well defined degree. In the case where the map is induced by
a map between the underlying spaces the degree behaves entirely
as we might expect:
\begin{theorem}
If $f\co\Sigma\rightarrow\Sigma'$ is a map between closed oriented
surfaces then
\[
\deg \Exp{k}{f} = (\deg f)^k.
\]
\label{surfacedegree.th}
\end{theorem}
Recall from~\cite{circles02} that the corresponding result for the 
circle is $\deg \Exp{k}{f} = 
\left(\deg f\right)^{\left\lfloor\frac{k+1}{2}\right\rfloor}$.
Both results are proved in the same way, by counting preimages of a 
generic point in the target with sign, and the difference lies in the 
fact that the circle is odd-dimensional. Orienting tangent spaces to
$\exps{k}\setminus\exps{k-1}$ as in~\eqref{surfaceorientation.eq}
requires the additional data of an order on the summands, and this
leads to cancellation among the preimages, whereas no such cancellation
occurs for surfaces.

As our last result for closed surfaces we construct a homotopy
model for the third finite subset space of a closed orientable
surface in section~\ref{model.sec}, and use this to 
completely calculate the integral cohomology groups of its third
finite subset space in section~\ref{exp3sigcohomology.sec}.
\begin{theorem}
\label{3rdorientable.th}
Let $\Sigma_g$ be a closed orientable surface of genus $g$. The cohomology
group $H^i(\expsig[g]{3})$ is given by 
\[
\begin{array}{|cr|ccc|} \hline
 & & \multicolumn{3}{c|}{$g$} \\
 &  &    0     &    1       &  \geq 2    \\ \hline
 &0 & \integer & \integer   & \integer  \\
 &1 &    0     &    0       &    0      \\
 &2 &    0     & \integer   & \integer^{\binom{2g}{2}} \\
i&3 &    0     & \integer^5 & \integer^{p(g)} \\
 &4 & \integer & \integer^4 & \integer^{q(g)}\oplus[\integer/2\integer]^{2g} \\
 &5 & \integer/2\integer & \integer^2\oplus\integer/2\integer
                            & \integer^{2g}\oplus\integer/2\integer\\
 &6 & \integer & \integer   & \integer  \\ \hline 
\end{array}
\]
in which 
\[
p(g) = \binom{2g}{3} + \binom{2g}{2} +\binom{2g}{1}  \qquad \mbox{and} \qquad
q(g) = \binom{2g}{2} + 1
\]
for $g\geq 2$. The Euler characteristic is
\[
\chi(\expsig{3}) = \frac{-4g^3 + 12g^2 -17g + 9}{3}.
\]
In particular $\chi(\Exp{3}{S^2})=3$, $\chi(\Exp{3}{T^2})=0$, and
$\chi(\Exp{3}{\Sigma_2})=-3$.
\end{theorem}
\noindent

\section{Lexicographic cell structures}
\label{lexicographic.sec}

We show that the finite subset spaces of connected finite $2$-complexes
admit cell structures based on the lexicographic ordering of $I^2$. 
After recalling some facts about finite subset spaces of graphs and 
establishing some conventions in section~\ref{Lconventions.sec}
we prove the cell structures exist in section~\ref{lexistence.sec}, and
calculate the boundary maps in section~\ref{boundary.sec}. We close in 
section~\ref{higherdimensions.sec} with a brief discussion on 
constructing lexicographic cell structures for the finite subset
spaces of higher dimensional complexes.

We remark that the existence of lexicographic cell structures is
chiefly of practical rather than theoretical importance, as the finite
subset spaces of a simplicial complex may be shown to have cell
structures using the machinery of simplicial sets 
(Jacob Mostovoy, private communication). Simplicial sets are 
described in May's book~\cite{may67} and Curtis's article~\cite{curtis71}.
Given a simplicial set $K$ we let \Exp{j}{K}\ be the simplicial
set whose $n$--simplices are subsets of size at most $j$ of the
$n$--simplices of $K$, and whose face and degeneracy operators
are the face and degeneracy operators of $K$ acting elementwise. 
Then if $X$ is the geometric realisation of $K$, \Exp{j}{X}\ will
be the geometric realisation of \Exp{j}{K}, showing that
triangulated spaces have triangulated finite subset spaces. The power of 
this method comes at the expense of difficulties with concrete
calculations: for example, the triangulations produced for
\Exp{3}{S^2}\
and \Exp{4}{S^2}\ have $77$ and $1039$ cells,\footnote{These numbers
are derived from cell counts sent by Jacob Mostovoy.} in contrast with 
the corresponding lexicographic cell structures which have $11$ and $23$.

Although their existence is not
of immediate theoretical interest, the   construction 
of lexicographic cell structures for the finite subset spaces
of connected $2$--complexes does
have an interesting theoretical consequence. Up to homotopy we may
assume that $X$ has a single vertex, and then the lexicographic
cell structure for \Exp{k}{X}\ is obtained
from that of \Exp{k-1}{X}\ by adding cells in dimensions $k-1\leq d\leq 2k$,
or $k\leq d\leq 2k$ if $X$ has no edges. Using a result of
Handel~\cite{handel00} this implies that \Exp{k}{X}\ is
$(k-2)$--connected, and $(k-1)$--connected if $X$ is simply connected.
In~\cite{vanish03} we use this to show that the same is true
without the dimension or finiteness restrictions.

\subsection{Conventions and definitions}
\label{Lconventions.sec}

In constructing cell structures for the finite subset spaces of a
connected finite $2$--complex $X$ we will work up to homotopy and
assume that $X$ has a single vertex, as we did with graphs in~\cite{graphs03}.
There we let
$\Gamma_n$ be the graph with a single vertex $v$ and $n$ edges
$e_1,\ldots,e_n$, and here we will identify the $1$--skeleton of $X$ with
$\Gamma_n$ for some $n\geq 0$ and denote its $2$--cells by
$f_1,\ldots,f_m$. Each cell $f_i$ has a
characteristic map $\phi_i$ from the unit disc $D^2\subseteq\complex$
to $X$ and we assume that $\phi_i$ sends $-1\in\partial D^2$ to $v$.

By Lemma~1 of~\cite{graphs03} \expg{k}{n}\ has a cell structure 
consisting of cells $\vcell{J}$, $\ecell{J}$ in which the indices are
 $n$--tuples
$J=(j_1,\ldots,j_n)$ of non-negative integers. The cells
\vcell{J}\ form a cell structure for \expgb{k}{n}, and the indexing 
$n$--tuple of the cell \vcell{J}\ or \ecell{J}\ containing the subset
$\Lambda$ in its interior is determined by the 
integers $j_i=|\Lambda\cap\intr e_i|$.
Letting $\tilde\Delta_j$ be the simplex
\[
\tilde\Delta_j = \{(x_1,\ldots,x_j) | 
                             0\leq x_1 \leq \cdots \leq x_j \leq 1 \},
\]
the domain of $\vcell{J}$, $\ecell{J}$ is the ball
$\tilde\Delta_J = \tilde\Delta_{j_1}\times\cdots\times\tilde\Delta_{j_n}$,
in which we omit any empty factors $\tilde\Delta_0$.

\begin{figure}[t]
\begin{center}
\leavevmode
\psfrag{kappa}{$\kappa$}
\psfrag{I2}{$I^2$}
\psfrag{D2}{$D^2$}
\includegraphics{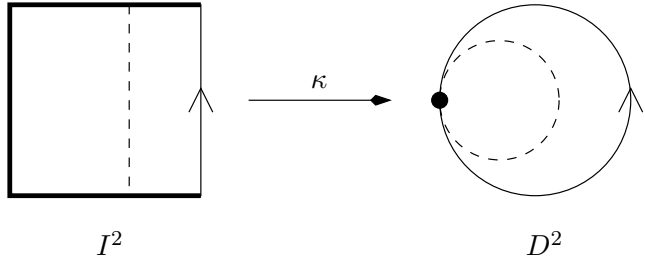}
\caption[The crunching map $\kappa$]{The crunching map 
        $\kappa\co I^2\rightarrow D^2$. The three sides of the
        square shown in bold map to the dot at $-1$, the unbolded
        edge maps to the boundary, and the dashed line goes to the
        dashed circle.}
\label{crunch.fig}
\end{center}
\end{figure}

The cell structures for \expg{k}{n}\ were constructed by using the
linear order on each edge to choose preferred lifts from \expg{j}{n}\
to $(\Gamma_n)^j$. To build cell structures for \Exp{k}{X}\ we
proceed similarly using the lexicographic order $\prec$, defined on
the unit square $I^2$ by
\begin{align*}
(x_1,y_1)&\prec(x_2,y_2) &\mbox{if $x_1 < x_2$, or if $x_1=x_2$ and $y_1<y_2$}.
\end{align*}
We transfer this order to $\intr D^2$ via a ``crunching'' map 
$\kappa\co I^2\rightarrow D^2$ which sends
the three sides $x=0$, $y=0$ and $y=1$ to 
$-1\in\partial D^2$ and maps the rest of the square homeomorphically onto 
the rest of the disc, preserving orientation. 
The exact form of $\kappa$ is unimportant, but for the sake of
concreteness we let
\[
\kappa(x,y) = x(1-e^{2\pi iy})-1.
\]
This sends the vertical segment $x=x_0$ to the circle with centre
$x_0-1$ and radius $x_0$, as shown in figure~\ref{crunch.fig}. 

As was the case with graphs our cells will be indexed by vectors
of integers corresponding to ordered partitions, but in contrast
with graphs the vectors arising from a $2$--cell will consist of 
positive integers only and their lengths will vary. 
We introduce some notation and terminology that will be
helpful in working with such partitions.

For each non-negative integer $m$ let 
$[m]=\{i\in\integer|1\leq i\leq m\}$. 
Given a vector $S$ of positive or non-negative integers we write
$\ell(S)=\ell$ if $S=(s_1,\ldots,s_\ell)$, and we define the norm of $S$
to be
\[
|S| = \sum_{i=1}^{\ell(S)} s_i.
\]
For each subset $\alpha$ of $[\ell(S)]$ we write 
$S|_\alpha$ for the $|\alpha|$--tuple obtained by restricting the index
set to $\alpha$.

The boundary of the cell corresponding to $S$ will
consist of two main contributions, one a sum of cells corresponding
to partitions obtained
by merging adjacent parts of $S$, and the second a sum of 
cells corresponding to partitions obtained by decreasing a part of $S$
by one. Accordingly we define
\[
\mu_i(S) = (s_1,\ldots,s_i+s_{i+1},\ldots,s_\ell)
\]
for $1\leq i\leq \ell(S)-1$, and 
\[
\drop_i(S) = (s_1,\ldots,s_i-1,\ldots,s_{\ell})
\]
for $1\leq i\leq \ell(S)$. Notice that $\mu_i$ decreases
the length of $S$ by one, while $\drop_i$ decreases the norm
of $S$ by one.

Finally, the length-decreasing contribution to the boundary will involve
the $(-1)$--binomial co-efficient $\sbinom{m}{r}$, which we recall
is a signed count of the  ways of choosing
$r$ elements from the set $[m]$.  Each choice is counted with the sign of the
permutation obtained by ordering the set $[m]$ so that the chosen
elements occur first in ascending order, followed by the unchosen
elements in ascending order.  The value of $\sbinom{m}{r}$ was
calculated in~section 3.2 of~\cite{graphs03} and is given by
\[
\sbinom{m}{r} = 
  \frac{1+(-1)^{r(m-r)}}{2}\binom{\lfloor m/2\rfloor}{\lfloor r/2\rfloor}.
\]

\subsection{Existence of lexicographic cell structures}
\label{lexistence.sec}

\begin{figure}[t]
\begin{center}
\leavevmode
\psfrag{1}{$1$}
\psfrag{2}{$2$}
\psfrag{3}{$3$}
\psfrag{(3)}{$(3)$}
\psfrag{(1,1,1)}{$(1,1,1)$}
\psfrag{(1,2)}{$(1,2)$}
\psfrag{(2,1)}{$(2,1)$}
\includegraphics{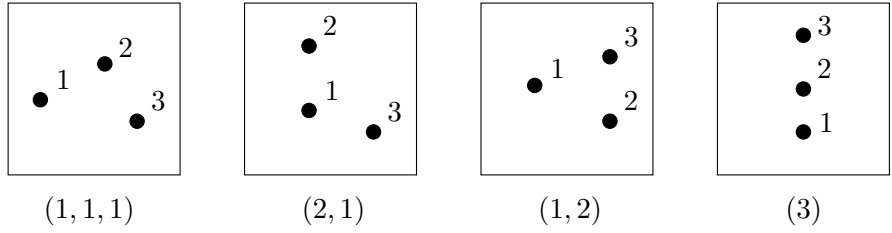}
\caption[Open cell decomposition of 
         $\Exp{3}{(\intr I^2)}\setminus\Exp{2}{(\intr I^2)}$]
        {Open cell decomposition of 
         $\Exp{3}{(\intr I^2)}\setminus\Exp{2}{(\intr I^2)}$.
         We decompose $\Exp{3}{(\intr I^2)}\setminus\Exp{2}{(\intr I^2)}$
         as a union of four open cells, corresponding to the four
         ordered partitions $(1,1,1)$, $(2,1)$, 
         $(1,2)$ and $(3)$ of three as a sum
         of positive integers.}
\label{lexico.fig}
\end{center}
\end{figure}

The first step in building a cell structure for $\Exp{k}{X}$ is to
form an open cell decomposition of 
$\Exp{j}{(\intr f_i)}\setminus\Exp{j-1}{(\intr f_i)}$, using the
fact that each $j$ element
subset of $\intr I^2$ has a unique lexicographically ordered
representative in $(\intr I^2)^j$. Figure~\ref{lexico.fig} illustrates
the idea for $k=3$. Generically, three points in $\intr I^2$ may be
ordered by their $x$ co-ordinates, as in the square labeled $(1,1,1)$,
and this gives an open $6$--ball
\[
\intr (\tilde\Delta_3\times \tilde\Delta_1\times 
            \tilde\Delta_1\times \tilde\Delta_1)
    =\{0<x_1<x_2<x_3<1\}\times\{0<y_1,y_2,y_3<1\}.
\]
If two of the $x$ co-ordinates are equal but the third is different there
are two possibilities, illustrated by the squares labeled $(2,1)$ and 
$(1,2)$ and corresponding to the open $5$-balls 
\begin{multline*}
\intr (\tilde\Delta_2\times\tilde\Delta_2\times\tilde\Delta_1) \\
       =\{0<x_1=x_2<x_3<1\}\times\{0<y_1<y_2<1\}\times\{0<y_3<1\},
\end{multline*}
\begin{multline*}
\intr (\tilde\Delta_2\times\tilde\Delta_1\times\tilde\Delta_2) \\
       =\{0<x_1<x_2=x_3<1\}\times\{0<y_1<1\}\times\{0<y_2<y_3<1\},
\end{multline*}
and when all three $x$ co-ordinates are equal as in the square labeled
$(3)$ we have the open $4$--ball
\[
\intr (\tilde\Delta_1\times\tilde\Delta_3) =
    \{0<x_1=x_2=x_3<1\}\times \{0<y_1<y_2<y_3<1\}.
\]
This gives a decomposition of 
$\Exp{3}{(\intr I^2)}\setminus\Exp{2}{(\intr I^2)}$ as a union of
four open cells, and more generally we obtain a decomposition of
$\Exp{j}{(\intr I^2)}\setminus\Exp{j-1}{(\intr I^2)}$ as a union 
of $2^{j-1}$ open cells, indexed by the ordered partitions
of $j$ as a sum of positive integers. A cell structure for \Exp{k}{X}\
is then obtained by taking products of such cells with the cells
$\vcell{J},\ecell{J}$ of \expg{k}{n}. 

Concretely, to each $j$ element subset 
$\Lambda=\{p_1,\ldots,p_j\}$ of $\intr I^2$ 
we associate the ordered partition 
$\S(\Lambda)=(s_1,\ldots,s_\ell)$ of $j$ arising from the equivalence
relation $p_i\sim p_j$ if $x_i=x_j$. This gives a partition
$\{\Lambda_1,\ldots,\Lambda_\ell\}$ of $\Lambda$ which may be ordered by
$\Lambda_a<\Lambda_b$ if $p\in\Lambda_a$, $p'\in\Lambda_b$, and $x<x'$, 
and we obtain $\S(\Lambda)$ by letting 
$s_i$ be the size of the $i$th part. Set
\[
\E(S) = \{\Lambda\in \Exp{|S|}{(\intr I^2)} | \S(\Lambda)=S\}
\]
for each vector $S$ of positive integers. It is clear that
the sets $\E(S)$ with $|S|=j$ form a partition of
$\Exp{j}{(\intr I^2)}\setminus\Exp{j-1}{(\intr I^2)}$, and
we claim further that 
$\E(S)$ is parameterised by the open $(|S|+\ell(S))$--ball 
$\intr(\tilde\Delta_S\times\tilde\Delta_{\ell(S)})$. 
To see this, note that if $\{p_1,\ldots,p_j\}\in\E(S)$ satisfies
$p_1\prec\cdots\prec p_j$ then
$(y_1,\ldots,y_j)$ is a point in 
$\intr (\tilde\Delta_{s_1}\times\cdots\times\tilde\Delta_{s_\ell})$ and 
the distinct
$x$ values $(x_{i_1},\ldots,x_{i_{\ell}})$ form a point in 
$\intr \tilde\Delta_{\ell(S)}$.  This leads to a map
\begin{equation}
\intr (\tilde\Delta_{s_1}\times\cdots\times\tilde\Delta_{s_\ell}
         \times\tilde\Delta_\ell) \rightarrow (\intr I^2)^j
          \rightarrow  \Exp{j}{(\intr I^2)}
\label{opencells.eq}
\end{equation}
hitting precisely $\E(S)$, and this map is injective since each 
$\Lambda\in\Exp{j}{(\intr I^2)}\setminus\Exp{j-1}{(\intr I^2)}$ has
a unique lexicographically ordered representative in $(\intr I^2)^j$.

Having found an open cell decomposition of $\Exp{j}{(\intr I^2)}$ we
now construct cells for $\Exp{k}{X}$. The map~\eqref{opencells.eq}
extends to a map 
$\tilde\Delta_S\times\tilde\Delta_\ell\rightarrow\Exp{j}{I^2}$
and we let $\fcell[i]{S}$ be the composition
\[
\tilde\Delta_S\times\tilde\Delta_{\ell(S)}\rightarrow
\Exp{|S|}{I^2}\xrightarrow{\Exp{|S|}{\kappa}}
\Exp{|S|}{D^2}\xrightarrow{\Exp{|S|}{\phi_i}}\Exp{|S|}{f_i},
\]
$\vfcell[i]{S}$ the composition
\[
\tilde\Delta_S\times\tilde\Delta_{\ell(S)}\xrightarrow{\fcell[i]{S}}
\Exp{|S|}{f_i}\xrightarrow{\cup\{v\}}\Exp{|S|+1}{(f_i,v)}.
\]
Since both
$\kappa|_{\intr I^2}$ and $\phi_i|_{\intr D^2}$ are homeomorphisms so
are $\Exp{k}{(\kappa|_{\intr I^2})}$ and $\Exp{k}{(\phi_i|_{\intr
D^2})}$, and it follows that the maps $\fcell[i]{S}$ are
injective on the interiors of their balls of definition 
and give an open cell decomposition of $\Exp{j}{(\intr f_i)}$. Further,
\vfcell[i]{S}\ is injective on 
$\intr(\tilde\Delta_S\times\tilde\Delta_{\ell(S)})$ as well,  since
$\cup\{v\}$ is injective on $\Exp{j}{(\intr f_i})$. We claim:

\begin{theorem}
Let $X$ be a connected finite $2$--complex with a single vertex $v$. Then
the maps
\begin{equation}
\bigl\{\vfcell[i_1]{S_1}\cup\cdots\cup\vfcell[i_p]{S_p}\cup\vcell{J}
   |i_1<\cdots<i_r\mbox{ and } |J|+{\textstyle\sum_q|S_{q}|} \leq k-1\bigr\}
\label{Xvcells.eq}
\end{equation}
are the characteristic maps of a cell structure for \Exp{k}{(X,v)}, and
these maps together with 
\begin{equation}
\bigl\{\fcell[i_1]{S_1}\cup\cdots\cup\fcell[i_p]{S_p}\cup\ecell{J}
   |i_1<\cdots<i_r\mbox{ and } |J|+{\textstyle\sum_q|S_{q}|} \leq k\bigr\}
\label{Xcells.eq}
\end{equation}
are the characteristic maps of a cell structure for $\Exp{k}{X}$. 
\label{lexicographic.th}
\end{theorem}

Notice that $\Exp{k}{X}$ is obtained from $\Exp{k-1}{X}$ by adding cells
in dimensions $k-1\leq d\leq 2k$, or $k\leq d\leq 2k$ if $X$ has no
edges. Handel~\cite{handel00} has shown that for path connected Hausdorff
$Y$ the inclusion map $\Exp{k}{Y}\hookrightarrow\Exp{2k+1}{Y}$ is zero
on all homotopy groups, and together these results imply that \Exp{k}{X}\
is $(k-2)$--connected, and $(k-1)$--connected if $X$ is simply connected. 
We show that this holds more generally for higher dimensional
complexes in~\cite{vanish03}.

\begin{proof}
Each map $\vfcell[i_1]{S_1}\cup\cdots\cup\vfcell[i_p]{S_p}\cup\vcell{J}$ 
in~\eqref{Xvcells.eq}
is a map from a 
$|J|+\sum_q|S_{q}|+\sum_q\ell(S_{q})$ dimensional ball
to \Exp{(|J|+\sum_q|S_{q}|+1)}{(X,v)}, and since 
$|J|+\sum_q|S_{q}|\leq k-1$ the
image of each map lies in $\Exp{k}{(X,v)}$. Similarly, each map
$\fcell[i_1]{S_1}\cup\cdots\cup\fcell[i_p]{S_p}\cup\ecell{J}$ occurring
in~\eqref{Xcells.eq} is a map from a ball of dimension 
$|J|+\sum_q|S_{q}|+\sum_q\ell(S_{q})$ with image lying in 
\Exp{k}{X}, and moreover the interior of this ball misses \Exp{k}{(X,v)}.
We show that each map is injective on the interior of its domain, that
the open cells partition \Exp{k}{X}, and that the boundary of each cell
lies on open cells of strictly smaller dimension.

To prove injectivity on interiors note that each basic cell 
\fcell[i_q]{S_{q}}, \ecell[i]{j_i}\ is injective on the interior of
its domain of definition, and that the images of these cells are disjoint
since they lie in disjoint sets $\Exp{|S_{q}|}{(\intr f_{i_q})}$, 
$\Exp{j_i}{(\intr e_i)}$. It follows that the restriction of the cupped
map $\fcell[i_1]{S_1}\cup\cdots\cup\fcell[i_p]{S_p}\cup\ecell{J}$ to
the interior of its ball of definition is injective also. Further, as
noted above this restriction misses \Exp{k}{(X,v)}, implying that
\[
\vfcell[i_1]{S_1}\cup\cdots\cup\vfcell[i_p]{S_p}\cup\vcell{J}
= (\cup\{v\})\circ
        (\fcell[i_1]{S_1}\cup\cdots\cup\fcell[i_p]{S_p}\cup\ecell{J})
\]
is injective on the interior of its ball of definition also. 

Each $\Lambda\in \Exp{k}{X}$ can be written uniquely as 
\[
\Lambda = (\Lambda\cap\{v\})\cup\bigcup_{i=1}^n(\Lambda\cap(\intr e_i))\cup
           \bigcup_{i=1}^m(\Lambda\cap(\intr f_i))
\]
and this data determines a unique open cell
from~\eqref{Xvcells.eq} or~\eqref{Xcells.eq} containing $\Lambda$. 
Let $i_1<\cdots<i_p$ be the indices $i$ such that 
$\Lambda\cap(\intr f_i)$ is non-empty, set 
$S_q=\S(\Lambda\cap(\intr f_{i_q}))$, and let $j_i=|\Lambda\cap(\intr e_i)|$.
Then the open 
cell $\vfcell[i_1]{S_1}\cup\cdots\cup\vfcell[i_p]{S_p}\cup\vcell{J}$
contains $\Lambda$ if $v\in\Lambda$, and the open cell
$\fcell[i_1]{S_1}\cup\cdots\cup\fcell[i_p]{S_p}\cup\ecell{J}$ contains
$\Lambda$ if $v\not\in\Lambda$. 
Moreover, this is the unique open cell containing
$\Lambda$: regarding 
$\vfcell[i_1]{S_1}\cup\cdots\cup\vfcell[i_p]{S_p}\cup\vcell{J}$ as
$\fcell[i_1]{S_1}\cup\cdots\cup\fcell[i_p]{S_p}\cup\ecell{J}\cup\{v\}$,
the map $\fcell[i]{S}$ can occur as a factor if and only if 
$\Lambda\cap(\intr f_i)$ is non-empty and $S=\S(\Lambda\cap(\intr f_i))$, 
\ecell{J}\ can occur as a factor if and only if 
$j_i=|\Lambda\cap(\intr e_i)|$ for each $i$, 
and $\{v\}$ can occur as a factor if
and only if $v\in\Lambda$. It follows that the open 
cells~\eqref{Xvcells.eq} and~\eqref{Xcells.eq} partition $\Exp{k}{X}$
as claimed.

It remains to check that taking boundaries decreases dimension. It
will suffice to do this for each factor cell, and since
we know the result for $\ecell{J},\vcell{J}$ by 
Lemma~1 of~\cite{graphs03} it will be
enough to check the cells \vfcell[i]{S}, \fcell[i]{S}. We do this in 
section~\ref{boundary.sec}, where we calculate 
$\partial\vfcell[i]{S}$ and $\partial\fcell[i]{S}$.
\end{proof}

\subsection{The boundary maps}
\label{boundary.sec}

In this section we indicate how to calculate the boundaries of the
cells~\eqref{Xvcells.eq} and~\eqref{Xcells.eq}, and complete the
proof of Theorem~\ref{lexicographic.th} by showing that taking
boundaries decreases dimension. Since 
\begin{align}
\partial(\sigma\cup\tau)&=\partial(\cup\circ(\sigma\times\tau)) \nonumber \\
             &= \cup_\sharp\partial(\sigma\times\tau) \nonumber \\
             &= \cup_\sharp(\partial\sigma\times\tau 
                         +(-1)^{\dim\sigma}\sigma\times\partial\tau)
                           \label{boundaryofproduct.eq}
\end{align}
we need only calculate the boundary of each basic cell occurring as
a factor and understand certain special cases of the cellular chain
map $\cup_\sharp$. The only nontrivial case that will arise in this
context is the case $\cup_\sharp(\ecell{J}\times\ecell{L})$ studied
in~\cite{graphs03}, so we will be able to confine our attention 
to calculating the boundaries of \vfcell[i]{S}\ and \fcell[i]{S}. 
To do this we calculate $\partial\fcell[i]{S}$ and obtain
$\partial\vfcell[i]{S}$ as $(\cup\{v\})_\sharp \partial\fcell[i]{S}$ by
removing all tildes.

The boundary of $\fcell[i]{S}$ will consist of three parts,
one that decreases the norm of $S$, a second that decreases the
length of $S$, and a third that moves points from the interior
of $f_i$ onto the boundary.
In preparation for stating the result we define each separately.
Recall from section~\ref{Lconventions.sec} that $\drop_a(S)$ is
the partition obtained by decreasing
the $a$th part of $S$ by one, and  that $\mu_a(S)$ is
the partition obtained by merging the $a$th and $(a+1)$th parts of $S$,
in other words
\begin{align*}
\drop_a(S) &= (s_1,\ldots,s_a -1, \ldots, s_\ell), \\
\mu_a(S) &= (s_1,\ldots,s_a+s_{a+1},\ldots,s_\ell).
\end{align*}
Since $\drop_a(S)$ has norm $|S|-1$ and length $\ell(S)$ the cell
$\fcell[i]{\drop_a(S)}$ has dimension one less than \fcell[i]{S},
and the same is true of $\fcell[i]{\mu_a(S)}$ since $|\mu_a(S)|=|S|$
and $\ell(\mu_a(S))=\ell(S)-1$. 
Let 
\begin{align*}
\dnu \vfcell[i]{S} &= -\sum_{a=1}^{\ell(S)} \frac{1+(-1)^{s_a}}{2} 
         (-1)^{\left| S|_{[a-1]}\right|}\vfcell[i]{\drop_a(s)}, \\
\dlambda\vfcell[i]{S} &= 
      \sum_{a=1}^{\ell(S)-1}(-1)^{(a-1)}\sbinom{s_a+s_{a+1}}{s_a}
                              \vfcell[i]{\mu_a(S)},
\end{align*}
and extend these maps to \fcell[i]{S}\ by
\begin{align*}
\dnu\fcell[i]{S} &= 2\dnu\vfcell[i]{S} -\widetilde{\dnu\vfcell[i]{S}}, 
\label{dgamma.eq} \\
\dlambda\fcell[i]{S} &= \widetilde{\dlambda\vfcell[i]{S}}.
\end{align*} 
The Greek letters $\nu$ and $\lambda$ are intended as mnemonics
for ``norm'' and ``length'' respectively, and as with graphs
a tilde over a chain $\tau$ means each cell $\vfcell[i]{S}$
in $\tau$ should be replaced with $\fcell[i]{S}$. 
The linear combinations $\dnu\vfcell[i]{S}$, $\dnu\fcell[i]{S}$
can be recognised as $\partial\vcell{S}$, $\partial\ecell{S}$ with
each $\sigma$ replaced by $\tau_i$, and it follows that $\dnu$
is a boundary operator. 
It is a consequence of the boundary calculation below with $X$ specialised to 
$S^2$ that $\dlambda$ is a boundary operator as well, and that these two 
operators commute.

For the third contribution to the boundary of $\fcell[i]{S}$ we
identify $\partial D^2$ with $\Gamma_1$, sending $v$ to $-1$ and
giving $\partial D^2$ the anti-clockwise orientation. Writing
$w_i\co\Gamma_1\rightarrow\Gamma_n$ for the attaching map of $f_i$,
there are chain maps $(\Exp{j}{w_i})_\sharp$ and we define
\begin{align*}
\dgamma\vfcell[i]{S} &= \vfcell[i]{S|_{[\ell-1]}}\cup
          (\Exp{s_{\ell}}{w_i})_\sharp\vcell{s_{\ell}}, \\
\dgamma\fcell[i]{S}  &= \widetilde{\dgamma\vfcell[i]{S}}.
\end{align*}
Then:
\begin{theorem}
The boundary of \fcell[i]{S}\ is given by
\[
\partial\fcell[i]{S} = \dnu\fcell[i]{S} + (-1)^{|S|}\dlambda\fcell[i]{S}
          +(-1)^{|S|+\ell-1}\dgamma\fcell[i]{S}.
\]
\label{boundarytau.th}
\end{theorem}

To calculate the boundary of a product cell
$\fcell[i_1]{S_1}\cup\cdots\cup\fcell[i_p]{S_p}\cup\ecell{J}$ 
use~\eqref{boundaryofproduct.eq} and observe that the only
nontrivial cases of $\cup_\sharp$ that occur come from the $\dgamma$ terms
and have the form
$\cup_\sharp(\fcell[i_1]{S_1}\times\cdots
\times\ecell{L}\times\cdots\times\fcell[i_p]{S_p}\times\ecell{J})$.
Since $\sigma\times\tau=\pm\tau\times\sigma$
we need only understand the case $\cup_\sharp(\ecell{L}\times\ecell{J})$, 
which is given by Section~3 of~\cite{graphs03}.

\begin{proof}[Proof of Theorem~\ref{boundarytau.th}]
The cell \fcell[i]{S}\ is a map from the ball 
$\tilde\Delta_{s_1}\times\cdots\times\tilde\Delta_{s_\ell}
\times\tilde\Delta_\ell$ and we
consider the effect of replacing an inequality with an equality
in a simplex in this product. There are two main cases, the simplices
$\tilde\Delta_{s_a}$ and $\tilde\Delta_\ell$, and each breaks into further
cases according to whether the equality occurs in the first, last, or 
an inner position. Four of the possibilities are
illustrated in figure~\ref{boundary.fig}.

\begin{figure}[t]
\begin{center}
\leavevmode
\psfrag{i}{(i)}
\psfrag{ii}{(ii)}
\psfrag{iii}{(iii)}
\psfrag{iv}{(iv)}
\includegraphics{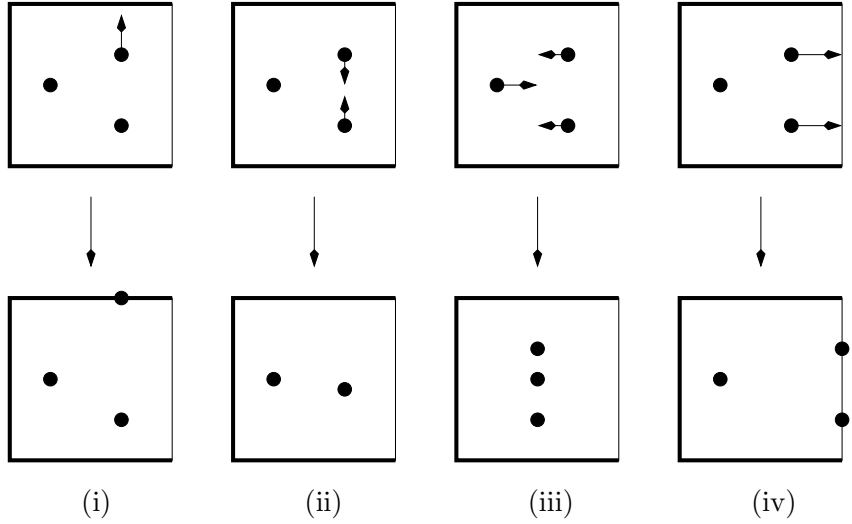}
\caption[Boundary calculations for finite subset spaces of $2$--complexes]
        {Some cells in the boundary of $\fcell{(1,2)}$. (i) Moving a
         point onto the top or bottom edge gives a set lying in 
         $\vfcell{(1,1)}$ (illustrated) or $\vfcell{(2)}$ (not shown). (ii)
         Moving two points with the same $x$ co-ordinate into co-incidence
         gives a point in $\fcell{(1,1)}$. (iii) Equating distinct $x$
         co-ordinates merges two columns, giving a point in $\fcell{(3)}$, 
         and (iv) moving the rightmost column onto the boundary gives
         a point in $\fcell{(1)}\cup\ecell{2}$.}
\label{boundary.fig}
\end{center}
\end{figure}

Replacing the first or last inequality in $\tilde\Delta_{s_a}$ sends a
point in the interior of $I^2$ onto the top or bottom edge of the
square, as shown
in figure~\ref{boundary.fig}(i). Since both of these edges map to 
$v$ in $X$ each case gives a point in $\vfcell[i]{\drop_a(S)}$, dropping
dimension by one, unless $s_a=1$ in which case it gives a point
in \vfcell[i]{S|_{[\ell]\setminus\{a\}}}, dropping
dimension by two. A generic point in $\vfcell[i]{\drop_a(S)}$ is
hit once by each face, and as was the case with graphs in~\cite{graphs03}
these contributions come with opposite signs if $s_a$ is odd and 
matching signs if $s_a$ is even, the sign negative if $a=1$.

Replacing an inner inequality in $\tilde\Delta_{s_a}$ merges two
points with the same $x$ co-ordinate, giving a point in 
$\fcell[i]{\drop_a(S)}$  as shown in figure~\ref{boundary.fig}(ii). 
A generic point in this cell is hit $s_a-1$ times, once for each
inequality, and as was the case with graphs these contributions
alternate in sign. They therefore cancel if $s_a$ is odd and leave
one over if $s_a$ is even, the leftover sign being positive if 
$a=1$. Putting these together with the contributions to 
$\vfcell[i]{\drop_a(S)}$
from the previous paragraph we get the $\dnu\fcell[i]{S}$ term in
$\partial\fcell[i]{S}$.

We now look at the contributions to the boundary from $\tilde\Delta_\ell$.
Replacing the first inequality with an equality moves the entire first column
of points onto the left edge of the square, which maps to $v$ in $X$. Thus
this face maps onto $\vfcell[i]{S|_{[\ell]\setminus\{1\}}}$, dropping
dimension by at least two
and so contributing nothing to $\partial\fcell[i]{S}$.

Replacing the $a$th inner inequality merges the $a$th and $(a+1)$th 
columns of points, giving a point in $\fcell[i]{\mu_a(S)}$ 
as shown in figure~\ref{boundary.fig}(iii). Generic points in 
this cell are each hit $\binom{s_a+s_{a+1}}{s_a}$ times, once
for each choice of $s_a$ points to come from the left,
and counting these with signs gives the
$(-1)$--binomial co-efficient $\sbinom{s_a+s_{a+1}}{s_a}$. Collecting 
these contributions together gives the $\dlambda\fcell[i]{S}$ term in
$\partial\fcell[i]{S}$.

Finally, we look at the contribution to the boundary from the last
face of $\tilde\Delta_\ell$, which moves the last column of points
onto the right edge of the square, as shown in
figure~\ref{boundary.fig}(iv). The points left in the interior of the
square lie in the cell $\fcell[i]{S|_{[\ell-1]}}$, which has dimension
$|S|+\ell-s_\ell-1$, and the points on the boundary lie in some cell
\vcell{J}\ or \ecell{J}\ of \expg{s_\ell}{n}, which has dimension
$s_\ell$. Thus this face maps onto cells of smaller dimensions
also. The signed sum of the dimension $s_\ell$ cells hit in
\expg{s_\ell}{n}\ is given by
$(\Exp{s_\ell}{w_i})_\sharp\ecell{s_\ell}$, and we get the
$\dgamma\fcell[i]{S}$ term of $\partial\fcell[i]{S}$.
\end{proof}

The boundary map $\dlambda$ can be simplified by choosing a suitable basis
for the chain groups over \rational, at the expense of making $\dnu$
mildly more complicated. Recalling from section~\ref{Lconventions.sec} that
the $(-1)$--binomial co-efficient $\sbinom{m}{r}$ is given by
\[
\sbinom{m}{r} = 
  \frac{1+(-1)^{r(m-r)}}{2}\binom{\lfloor m/2\rfloor}{\lfloor r/2\rfloor}
\]
we let
\[
\qvfcell{S} = (\lfloor s_1/2\rfloor!\cdots\lfloor s_\ell/2\rfloor!)
                   \vfcell{S}.
\]
Then 
\begin{align*}
{\textstyle\sbinom{s_a+s_{a+1}}{s_a}} \qvfcell{\mu_a(S)} 
   &=\frac{1+(-1)^{s_a s_{a+1}}}{2}\,
     \frac{\lfloor (s_a+s_{a+1})/2\rfloor!}
          {\lfloor s_a/2\rfloor!\lfloor s_{a+1}/2\rfloor!}\,
          (\lfloor s_1/2\rfloor!\cdots\lfloor s_\ell/2\rfloor!)
              \vfcell{\mu_a(S)}\\
   &=\frac{1+(-1)^{s_a s_{a+1}}}{2}\qvfcell{\mu_a(S)},
\end{align*}
so that with respect to this basis we have
\begin{equation}
\label{rationalboundary.eq}
\dlambda\qvfcell{S} 
    = \sum_{a=1}^{\ell(S)-1}\frac{1+(-1)^{s_a s_{a+1}}}{2}
                        (-1)^{(a-1)} \qvfcell{\mu_a(S)}.
\end{equation}
This simplification comes at the small price of adding a factor
of $\lfloor s_a/2\rfloor$ to the $a$th term of $\dnu$ when $s_a$
is even, giving
\[
\dnu \qvfcell{S} = -\sum_{a=1}^{\ell(S)} \frac{1+(-1)^{s_a}}{4} 
   (-1)^{\left| S|_{[a-1]}\right|} s_a \qvfcell{\drop_a(s)}.
\]
This mild complication in $\dnu$ will not concern us, since our main
use of~\eqref{rationalboundary.eq} will be in calculating the
rational homology of $\expstwob{k}/\expstwob{k-1}$, where the boundary
map consists only of $\dlambda$. We note further that $\qvfcell{S}=\vfcell{S}$
if each entry of $S$ is at most three. 

\subsection{Higher dimensions}
\label{higherdimensions.sec}

We make a brief digression on constructing lexicographic cell
structures for the finite subset spaces of higher dimensional
complexes. As was the case with $2$--complexes, the lexicographic
ordering on $I^n$ may be used to construct an open cell decomposition
of $\Exp{j}{(\intr I^n)}$, and we might hope to use this to build a
cell structure for $\Exp{k}{X}$ by taking products over cells of $X$ 
as before. However, ensuring that boundaries
decrease dimension under this scheme 
appears to require compatibility between the
ordering on the interior of a cell and the orderings on cells in its
boundary. This condition comes for free in the case of $2$--complexes
or wedges of spheres but appears to require orchestration in general.

To each lexicographically ordered 
sequence $p_1,\ldots,p_j$ of $j$ points in the
interior of $I^n$ we associate a $(j-1)$--tuple of integers
by setting $j_i=m-1$ if $p_i$ and $p_{i+1}$ are
distinguished by their $m$th co-ordinate but not their $(m-1)$th. This
$(j-1)$--tuple describes the dependencies among the co-ordinates of the
$p_i$ and as before we consider the collection of lexicographically
ordered subsets corresponding to a fixed tuple $J$. This is easily
seen to be parameterised by an open ball of dimension $nj-|J|$,
a product of open simplices of the form $\intr\tilde\Delta_a$,
and this gives an open cell decomposition of $\Exp{j}{(\intr I^n)}$. 
We again transfer this to $\Exp{j}{(\intr B^n)}$ by collapsing 
the entire boundary to a point, with the exception of the
sole face $x=1$; the purpose of this is to avoid constraints between points
in the interior and points on the boundary which would otherwise arise
along faces of the form shown in figure~\ref{boundary.fig}(i).

To illustrate the compatibility requirement let $n=3$ and give the face $x=1$ 
the reverse lexicographic ordering
\begin{align*}
(1,y_1,z_1)&\prec(1,y_2,z_2) 
     &\mbox{if $z_1 < z_2$, or if $z_1=z_2$ and $y_1<y_2$}.
\end{align*}
Three element subsets of $\intr I^3$ corresponding to the $2$--tuple $(2,2)$ 
satisfy $x_1=x_2=x_3$, $y_1=y_2=y_3$
and $z_1<z_2<z_3$, and are parameterised by the open $5$--ball
$\intr(I\times I\times\tilde\Delta_3)$. Pushing them onto the boundary
we obtain subsets corresponding to the partition
$(1,1,1)$, which is the $2$--tuple $(0,0)$ in the labeling used
here. This is parameterised by the open $6$--ball
$\intr (I^3\times\tilde\Delta_3)$, showing that 
dimensions can jump if the ordering on the boundary does not agree with
that coming from the interior.

\section{The two-sphere}
\label{twosphere.sec}

\subsection{Introduction}

We now use the lexicographic cell structures developed in the
previous section to study the finite subset spaces of the two-sphere.
We use the standard cell structure for \stwo\ with a single
two-cell, and to avoid clutter we will simply write $S$ and $\tilde S$ for
\vfcell{S}\ and \fcell{S}\ respectively. Where the partition is
written out explicitly we will use round brackets for $\vfcell{S}$ and
square brackets for $\fcell{S}$.

We begin in section~\ref{spherechains.sec}
by describing the chain complexes of \expstwob{k}\ and \expstwo{k}
and justifying the assertion of section~\ref{boundary.sec}
that $\dlambda$ and $\dnu$ are commuting boundary
operators. 
We look at the finite subset spaces of $\stwo$ for several small
values of $k$ in section~\ref{smallk.sec}, and then prove
Theorems~\ref{s2rational.th} and~\ref{s2integral.th} in 
section~\ref{spherehomology.sec}.
The arguments bear a striking resemblance to those used to calculate
the homology of \expg{k}{n}\ in~\cite{graphs03}. 
The based space \expstwob{k}\ again
plays a prominent role, being easier to understand and having
cleaner results, and we make use of
the exact sequences of the pairs $(\expstwob{k},\expstwob{k-1})$ and
$(\expstwo{k},\expstwob{k})$. The main step is the calculation
of $H_*(\expstwob{k}/\expstwob{k-1};\rational)$. This involves the 
study of a combinatorially defined 
finite complex which may be recognised as
the $(k-2)$--cube complex of~\cite{graphs03} with some signs changed and 
certain edges ``clipped''.

\subsection{The chain complexes}
\label{spherechains.sec}

\begin{figure}[t]
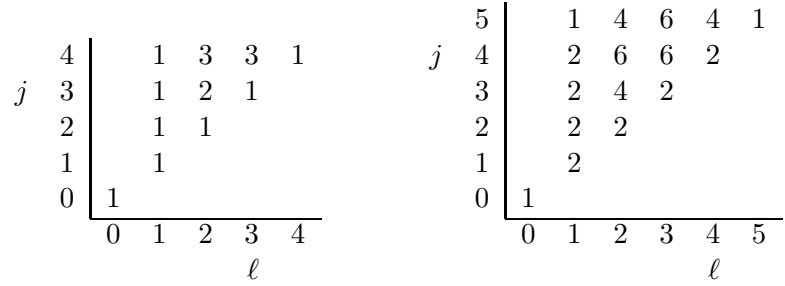

\leavevmode
\begin{center}
\begin{tabular}{ccc}
\begin{tabular}{rr|rrrrr}
\multicolumn{7}{r}{}               \\
       & 4 &   & 1 & 3 & 3 & 1   \\ 
$j$    & 3 &   & 1 & 2 & 1 &     \\
       & 2 &   & 1 & 1 &   &     \\
       & 1 &   & 1 &   &   &     \\
       & 0 & 1 &   &   &   &     \\ \cline{3-7}
\multicolumn{3}{r}
             0 & 1 & 2 & 3 & 4   \\
\multicolumn{6}{r}     {$\ell$} & 
\end{tabular} & $\quad$ &
\begin{tabular}{rr|rrrrrr}
       & 5 &   & 1 & 4 & 6 & 4 & 1 \\
$j$    & 4 &   & 2 & 6 & 6 & 2 &   \\ 
       & 3 &   & 2 & 4 & 2 &   &   \\
       & 2 &   & 2 & 2 &   &   &   \\
       & 1 &   & 2 &   &   &   &   \\
       & 0 & 1 &   &   &   &   &   \\ \cline{3-8}
\multicolumn{3}{r}
             0 & 1 & 2 & 3 & 4 & 5 \\
\multicolumn{7}{r}     {$\ell$} & 
\end{tabular}
\end{tabular}
\caption[Cell counts for $\expstwob{5}$ and $\expstwo{5}$]
        {Cell counts for $\expstwob{5}$ and $\expstwo{5}$. 
         Left: Cells of \expstwob{k}\
         may be arranged into a double complex with $\binom{j-1}{\ell-1}$
         cells in position $(\ell,j)$, $1\leq \ell\leq j\leq k-1$,
         and one in position $(0,0)$. Since the alternating sum along
         a row of Pascal's triangle is zero we see immediately that
         $\chi(\expstwob{k})=2$ for $k\geq 2$. Right: Adding cells to get
         \expstwo{k}\ doubles all rows other than the bottom one
         and adds an extra row of Pascal's triangle on top, giving
         $\chi(\expstwo{k})=3$ for $k\geq 2$.}
\label{cellcounts.fig}
\end{center}
\end{figure}

By Theorem~\ref{lexicographic.th} \expstwob{k}\ has a cell structure
with a single vertex in which the remaining cells are 
in one-to-one correspondence
with the ordered partitions of the positive integers $j\leq k-1$,
and \expstwo{k}\ has a cell structure 
obtained by adding additional cells corresponding
to the ordered partitions of $j\leq k$.  
Let $\S_{\ell,j}$ and $\tilde\S_{\ell,j}$ be two copies of 
the free abelian group of rank $\binom{j-1}{\ell-1}$
generated by the ordered partitions of $j$ of length $\ell$. Each 
generator corresponds to a cell of dimension $j+\ell$ and we
note that $\dnu$ and $\dlambda$ map $\S_{\ell,j}$, $\tilde\S_{\ell,j}$
to $\S_{\ell,j-1}$, $\S_{\ell,j-1}\oplus\tilde\S_{\ell,j-1}$ and 
$\S_{\ell-1,j}$, $\S_{\ell-1,j}\oplus\tilde\S_{\ell-1,j}$ respectively. 
Since the boundary
map is given by $\partial = \dnu + (-1)^j\dlambda$ it follows that
the cellular chain groups of \expstwob{k} and \expstwo{k} may be arranged
into double complexes with cell counts as shown in 
figure~\ref{cellcounts.fig}, and that $\dnu,\dlambda$ are commuting
boundary operators as claimed. Moreover, we see immediately that
for $k\geq 2$ we have $\chi(\expstwob{k})=2$ and $\chi(\expstwo{k})=3$.

\begin{figure}[t]
\leavevmode
\begin{center}
\psfrag{(1,1,1,1)}{$(1,1,1,1)$}
\psfrag{(2,1,1)}{$(2,1,1)$}
\psfrag{(1,2,1)}{$(1,2,1)$}
\psfrag{(1,1,2)}{$(1,1,2)$}
\psfrag{(1,3)}{$(1,3)$}
\psfrag{(2,2)}{$(2,2)$}
\psfrag{(3,1)}{$(3,1)$}
\psfrag{(4)}{$(4)$}
\includegraphics{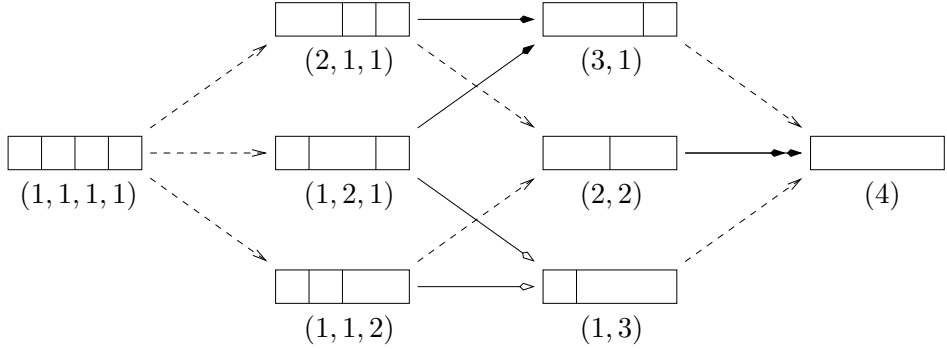}
\caption[The clipped cube complex]
        {The clipped cube complex; compare with figure~2 
         of~\protect\cite{graphs03},
         regarding the vertical dividing bars as the set $\{1,2,3\}$.
         The lattice of partitions of $4$ forms a $3$--dimensional
         cube, and the boundary of $S$ is a signed sum of its 
         neighbours with fewer parts, omitting those obtained by
         amalgamating adjacent parts of odd size. 
         Positive terms are indicated by solid arrowheads and negative
         terms by empty arrowheads, and the omitted terms (the ``clipped''
         edges) are shown by dotted arrows. The double-headed arrow from
         $(2,2)$ to $(4)$ shows where $\dq$ differs from $\dlambda$, as
         $\dlambda(2,2)=2(4)$ but $\dq(2,2)=(4)$. We see that the homology
         of the  clipped
         cube complex is $0$ in dimensions $1$ and $2$ and  
         \integer\ in dimensions $3$ and $4$,
         generated by $(2,1,1)-(1,2,1)+(1,1,2)$ and $(1,1,1,1)$ respectively,
         and that $\expstwob{5}/\expstwob{4}$ has an additional 
         $\integer/2\integer$ summand generated by $(4)$.          }
\label{clippedcube.fig}
\end{center}
\end{figure}

In working with these complexes we will focus mainly on the chain 
complexes $(\S_{*,j},\dlambda)$ consisting of a single 
row, which amounts to working with the quotient spaces 
$\expstwob{j+1}/\expstwob{j}$. The $(-1)$--binomial 
co-efficients appearing in $\dlambda$ will be extraneous to our purpose so
for clarity we will use the simplified boundary 
map~\eqref{rationalboundary.eq}, which we recall comes from using
bases for $\S_{\ell,j}$ 
obtained by individually scaling each cell. Abusing notation
we write~\eqref{rationalboundary.eq} as
\[
\dq S = \sum_{a=1}^{\ell}\frac{1+(-1)^{s_a s_{a+1}}}{2}
                        (-1)^{(a-1)} \mu_a(S).
\]
By analogy with the cube complex of section~2.3 of~\cite{graphs03} we call
the complex $(\S_{*,j},\dq)$ the ``clipped cube complex'', as the 
lattice of partitions of $j$, graded by length and 
partially ordered by refinement, forms
a ${(j-1)}$--dimensional cube, and the boundary of $S$ 
is again a linear combination
of its neighbours of smaller degree. 
The adjective comes
from the fact that $\dq S$ omits the partitions obtained from $S$
by amalgamating adjacent parts of odd size, so that these edges may
be regarded as removed or ``clipped''.  This is illustrated in
figure~\ref{clippedcube.fig}, which should be compared with
figure~2 of~\cite[p. 885]{graphs03}.
We note that $\dlambda$ and $\dq$ co-incide 
if each entry of $S$ and $\mu_a(S)$
have size at most three, which holds for the top two chain groups of
the complex. We will use this fact to calculate the top integral
homology groups of $\expstwob{k}$ and \expstwo{k}.
The distinction between $\dlambda$ and $\dq$ occurs for the first
time in the three dimensional clipped cube complex and can be seen in
figure~\ref{clippedcube.fig}.

An understanding of the clipped cube complex allows us to 
calculate the rational homology of $\expstwob{k}$ inductively, 
using  the long exact sequence of the pair 
$(\expstwob{k},\expstwob{k-1})$, 
and we pass
to an understanding of $\expstwo{k}$ using the
long exact sequence of the pair 
$(\expstwo{k},\expstwob{k})$. Our tool in doing so 
is the observation that
the natural correspondence $\tilde{S}\leftrightarrow S$ 
between the tilded cells of \expstwo{k}\ 
and the cells of \expstwob{k+1}\ has the following algebraic 
consequence:
\begin{lemma}
The spaces \expstwob{k+1} and $\expstwo{k}/\expstwob{k}$ have
the same homology.
\label{samehomology.lem}
\end{lemma}
\begin{proof} 
Comparing \expstwob{k+1} with  the quotient space
$\expstwo{k}/\expstwob{k}$ we see that one has a cell structure 
consisting of a single vertex, cells $S$ satisfying $|S|\leq k$, and
boundary map
\[
\partial S = {\dnu S} + (-1)^{|S|}{\dlambda S},
\]
and the second has a cell structure consisting
of a single vertex, cells $\tilde S$ satisfying
$|\tilde{S}|\leq k$, and boundary map
\[
\partial \tilde{S} = -\widetilde{\dnu S} + (-1)^{|S|}\widetilde{\dlambda S}.
\]
To account for the minus sign in the first term of the second map we
twist the natural correspondence $\tilde{S}\leftrightarrow S$ by
$(-1)^{|S|}$, obtaining an (algebraically defined) map 
$S\mapsto (-1)^{|S|}\tilde{S}$ inducing an isomorphism of 
chain complexes.
\end{proof}

We emphasise the fact that we have proved this result algebraically,
not topologically. The set map 
\[
\expstwob{k+1}\setminus\bigl\{\{*\}\bigr\}
\rightarrow \expstwo{k}/\expstwob{k} :
\Lambda\mapsto\Lambda\setminus\{\ast\}
\]
is discontinuous at each $\Lambda$ of size $k$ or less, and the map
\[
\expstwo{k}\rightarrow\expstwob{k+1} :
\Lambda\mapsto\Lambda\cup\{\ast\}
\]
does not descend to a continuous map on the quotient
$\expstwo{k}/\expstwob{k}$.

\subsection{Small values of $k$}
\label{smallk.sec}

In this section we calculate the homology of $\expstwob{k}$ and
$\expstwo{k}$ for several small
values of $k$. Our aim is both to build familiarity with the chain complexes
and to show by example that Theorem~\ref{s2integral.th} does not give the
full story on the integral homology of $\expstwo{k}$.

Since $\expstwo{1}$ and $\expstwob{2}$ are both simply $\stwo$ with
its standard cell structure we begin with \expstwo{2} and \expstwob{3}. 
The lexicographic cell structure for  \expstwo{2}\ 
has one cell $[1,1]$ in dimension
four, one cell $[2]$ in dimension three, two cells $[1]$ and $(1)$ 
in dimension two, no cells in dimension one and a single vertex. The
chain complex is
\[
0\longrightarrow \integer \decarrow{0} \integer 
\xrightarrow{\qbinom{1}{-2}}
\integer\oplus\integer \longrightarrow 0 \longrightarrow \integer 
\longrightarrow 0,
\]
and the homology is clearly
\[
H_i(\expstwo{2}) =\begin{cases}
                  \integer & i=0,2,4, \\
                  0        & \mbox{else}.
                  \end{cases}
\]
Mapping \expstwo{2}\ to \expstwob{3}\ by adding the vertex $*$ to 
each set rounds the brackets on the cells $[1,1]$, $[2]$ and $[1]$,
giving the chain complex
\[
0\longrightarrow \integer \decarrow{0} \integer 
\decarrow{-1}
\integer \longrightarrow 0 \longrightarrow \integer 
\longrightarrow 0
\]
with homology
\begin{equation}
H_i(\expstwob{3}) =\begin{cases}
                  \integer & i=0,4, \\
                  0        & \mbox{else}.
                  \end{cases}
\label{s4homology.eq}
\end{equation}
Geometrically, $\expstwo{2}=\sym{2}{\stwo}=\sym{2}{\complex P^1}$ 
may be identified with
$\complex P^2$ via the map
\[
\bigl\{[z_0,w_0],[z_1,w_1]\bigr\}
             \mapsto [z_0 z_1, -(z_0 w_1 + z_1 w_0), w_0 w_1],
\]
where the co-ordinates of the right-hand side are the co-efficients $[a,b,c]$
of the homogeneous quadratic
$(z_0 W - w_0 Z)(z_1 W - w_1 Z)$
vanishing at $[z_0,w_0]$ and $[z_1,w_1]$. In this picture 
the quotient map
$\expstwo{2}\xrightarrow{\cup\{*\}}\expstwob{3}$ 
induces the equivalence relation 
$[z^2,-2zw,w^2]\sim [z,-w,0]$, identifying the degree two curve
$b^2-4ac=0$ with the degree one curve $a+b=0$.  As seen
in~\eqref{s4homology.eq} this kills the second homology and the
resulting space $\expstwob{3}$ is homotopy equivalent to $S^4$.

To calculate the homology of \expstwob{4}\ and \expstwo{3}\ we start with
the clipped cube complex
\[
\begin{CD}
(1,1,1) @>0>> (2,1) \\
@VV0V        @VV1\quad .V  \\
(1,2)   @>1>> (3)    
\end{CD}
\]
This has homology groups $\integer$ in the top and middle dimensions, 
generated by $(1,1,1)$ and $(2,1)-(1,2)$ respectively, so
\[
\tilde{H}_i(\expstwob{4}/\expstwob{3}) = 
  \begin{cases}
         \integer & i=5,6, \\
          0       &\mbox{else}.
  \end{cases}
\]
The long exact sequence of the pair $(\expstwob{4},\expstwob{3})$ then gives
\[
0\longrightarrow  \tH_6(\expstwob{4}) \longrightarrow 
       \tH_6(\expstwob{4}/\expstwob{3})  \longrightarrow 0 
\]
and
\begin{equation}
0\longrightarrow  \tH_5(\expstwob{4}) \longrightarrow 
       \integer \decarrow{\partial}
       \integer  \longrightarrow 
       \tH_4(\expstwob{4})\longrightarrow 0,
\label{LESof(4,3).eq}
\end{equation}
the remaining segments having the form 
$0\longrightarrow \tilde{H}_i(\expstwob{4})\longrightarrow 0$. The
boundary map 
$\partial\co \tH_5(\expstwob{4}/\expstwob{3})\rightarrow \tH_4(\expstwob{3})$
is $\dnu$ which sends $(2,1)-(1,2)$ to $-2(1,1)$, so the 
middle map in~\eqref{LESof(4,3).eq} is multiplication by $-2$. Hence
\[
\tH_i(\expstwob{4}) =\begin{cases}
                  \integer           &   i= 6, \\
                  \integer/2\integer &   i=4, \\
                  0        & \mbox{else},
                  \end{cases}
\]
in which the nontrivial groups are generated by $(1,1,1)$ and $(1,1)$
respectively. 

We now use Lemma~\ref{samehomology.lem} and the exact sequence
of the pair $(\expstwo{3},\expstwob{3})$ to calculate the homology 
of \expstwo{3}. The segments where $\tH_i(\expstwo{3})$ is not
bracketed by zeroes are
\[
0\longrightarrow  \tH_6(\expstwo{3}) \longrightarrow 
       \tH_6(\expstwo{3}/\expstwob{3})  \longrightarrow 0, 
\]
giving $\tH_6(\expstwo{3})\cong\integer$ generated by $[1,1,1]$, and
\[
0\longrightarrow \integer \longrightarrow \tH_4(\expstwo{3})
 \longrightarrow \integer/2\integer \longrightarrow 0.
\]
To see that this second sequence splits consider the cycle
$[1,1]-2(1,1)$. Under the right arrow this maps to the generator $[1,1]$
of $\tH_4(\expstwo{3}/\expstwob{3})$, and the fact that the boundary
of $[2,1]-[1,2]$ is $2([1,1]-2(1,1))$  shows it has order two.
Hence
\[
\tH_i(\expstwo{3}) =\begin{cases}
                  \integer           &   i= 6, \\
                  \integer\oplus\integer/2\integer &   i=4, \\
                  0        & \mbox{else},
                  \end{cases}
\]
with generators $[1,1,1]$ in dimension six and $[1,1]$ and $(1,1)$ in
dimension four. Using the universal co-efficient
theorem we see that this result agrees
with the $g=0$ case of Theorem~\ref{3rdorientable.th}. 

From figure~\ref{clippedcube.fig} the homology of 
$\expstwob{5}/\expstwob{4}$ is $\integer$ in dimensions
seven and eight and $\integer/2\integer$ in dimension five.
The long exact sequence of $(\expstwob{5},\expstwob{4})$ has
two nontrivial segments
\[
0\longrightarrow  \tH_8(\expstwob{5}) \longrightarrow 
       \tH_8(\expstwob{5}/\expstwob{4})  \longrightarrow 0
\]
and
\[
0\longrightarrow  \tH_7(\expstwob{5}) \longrightarrow 
       \integer \decarrow{\partial}
       \integer  \longrightarrow 
       \tH_6(\expstwob{5})\longrightarrow 0
\]
as above, and an additional nontrivial segment
\[
0\longrightarrow\tH_5(\expstwob{5})\longrightarrow \integer/2\integer
 \decarrow{\partial}\integer/2\integer \longrightarrow
 \tH_4(\expstwob{5})\longrightarrow 0.
\]
The first boundary map sends the generator 
$(2,1,1)-(1,2,1)+(1,1,2)$ to $-3(1,1,1)$, inducing multiplication 
by $-3$, and the second sends the generator $(4)$ to $-(3)$. 
Since $\partial (2,1)=-(1,1)-(3)$ the second map is an 
isomorphism and
\[
\tH_i(\expstwob{5}) =\begin{cases}
                  \integer           &   i= 8, \\
                  \integer/3\integer &   i=6, \\
                  0        & \mbox{else}.
                  \end{cases}
\]
Again using Lemma~\ref{samehomology.lem} and 
the exact sequence of $(\expstwo{4},\expstwob{4})$
we get
\begin{gather}
0\longrightarrow  \tH_8(\expstwo{4}) \longrightarrow 
       \tH_8(\expstwo{4}/\expstwob{4})  \longrightarrow 0, \nonumber \\
0\longrightarrow \integer \longrightarrow \tH_6(\expstwo{4})
 \longrightarrow \integer/3\integer \longrightarrow 0,
\label{Z+Z/3Zseq.eq}
\end{gather}
and
\[
0\longrightarrow \integer/2\integer \longrightarrow 
 \tH_4(\expstwo{4})\longrightarrow 0 .
\]
The $6$--cycle $[1,1,1]-2(1,1,1)$ maps to the generator of the
$\integer/3\integer$ term and three times it is the boundary 
of $[2,1,1]-[1,2,1]+[1,1,2]$, so the second sequence splits
and
\[
\tH_i(\expstwo{4}) =\begin{cases}
                  \integer           &   i= 8, \\
                  \integer\oplus\integer/3\integer &   i=6, \\
                  \integer/2\integer       & i=4, \\
                  0                        & \mbox{else}.
                  \end{cases}
\]

We calculate one more example to show that the torsion becomes
increasingly complicated as $k$ increases. The example will 
require understanding the four dimensional clipped cube complex,
shown in figure~\ref{4Dclippedcube.fig},
and will further illustrate the ideas used in the general case.

\begin{figure}
\leavevmode
\begin{center}
\[
\begin{array}{ccc}
(1,1,1,1,1) \quad & 
\begin{aligned}
                 && (3,1,1) \quad\; &&                  \\
 (2,1,1,1) \quad && (1,3,1) \quad\; && (4,1) \quad      \\
 (1,2,1,1) \quad && (1,1,3) \quad\; && (1,4) \quad      \\
 (1,1,2,1) \quad && (2,2,1) \quad\; && (3,2) \quad      \\
 (1,1,1,2) \quad && (1,2,2) \quad\; && (2,3) \quad      \\
                 && (2,1,2) \quad\; &&       
\end{aligned} &
\; (5)
\end{array}
\]
\caption[Generators of the four dimensional clipped cube complex]
        {Generators of the four dimensional clipped cube complex.}
\label{4Dclippedcube.fig}
\end{center}
\end{figure}

As always, the boundary
of a generator with all entries odd is zero, so we may
consider the restriction of the boundary $\dlambda$ to the span of
the generators with one or more
even entries. Moreover, the boundaries of $(2,1,1,1)$ and its
permutations lie in the span of $(3,1,1)$ and its permutations,
so we may regard $\dlambda_4$ as a map to this subspace. With these
conventions and the bases ordered as shown the matrices of the
boundary maps are
\begin{align*}
\dlambda_4 &= 
\begin{bmatrix}
1 &  1 &  0 &  0 \\
0 & -1 & -1 &  0 \\
0 &  0 &  1 &  1 
\end{bmatrix}, &
\dlambda_3 &=
\begin{bmatrix}
 2 &  0 &  0 \\
 0 & -2 &  0 \\
 0 &  1 &  1 \\
-1 &  0 & -1 
\end{bmatrix}, &
\dlambda_2 & =
\begin{bmatrix}
1 \\ 1 \\ 2 \\ 2
\end{bmatrix}.
\end{align*}
Clearly, both $\dlambda_4$ and $\dlambda_2$ are surjective 
and $\dlambda_3$ is injective. 
The kernel of $\dlambda_2$ is equal to $\{2x+2y+z+w=0\}$ and by choosing
a suitable basis for this subspace it can be seen that the image of
$\dlambda_3$ has index two. In the now familiar pattern the long exact 
sequence of the pair $(\expstwob{6},\expstwob{5})$ gives 
\begin{gather*}
0\longrightarrow  \tH_{10}(\expstwob{6}) \longrightarrow 
       \tH_{10}(\expstwob{6}/\expstwob{5})  \longrightarrow 0, \\
0\longrightarrow  \tH_9(\expstwob{6}) \longrightarrow 
       \integer \decarrow{\partial}
       \integer  \longrightarrow 
       \tH_8(\expstwob{6})\longrightarrow 0,
\end{gather*}
and an additional segment
\[
0\longrightarrow \tH_7(\expstwob{6}) \longrightarrow
 \integer/2\integer \decarrow{\partial}\integer/3\integer
 \longrightarrow \tH_6(\expstwob{6})\longrightarrow 0.
\]
The boundary map in the third sequence is necessarily zero, and
the boundary map in the second sequence takes the generator
$(2,1,1,1)-(1,2,1,1)+(1,1,2,1)-(1,1,1,2)$ to $-4(1,1,1,1)$.
Hence
\[
\tH_i(\expstwob{6}) =\begin{cases}
                  \integer           &   i= 10, \\
                  \integer/4\integer &   i=8, \\
                  \integer/2\integer &   i=7, \\
                  \integer/3\integer &   i=6, \\
                  0        & \mbox{else}.
                  \end{cases}
\]

Finally, from the long exact sequence of the pair
$(\expstwo{5},\expstwob{5})$ we extract $\tH_{10}(\expstwo{5})\cong\integer$,
$\tH_7(\expstwo{5})\cong\integer/2\integer$, and two short exact sequences
\begin{gather*}
0\longrightarrow \integer \longrightarrow \tH_8(\expstwo{5})
 \longrightarrow \integer/4\integer \longrightarrow 0, \\
0\longrightarrow \integer/3\integer \longrightarrow \tH_6(\expstwo{5})
 \longrightarrow \integer/3\integer \longrightarrow 0.
\end{gather*}
As usual the first is split by $[1,1,1,1]-2(1,1,1,1)$, and 
as in~\eqref{Z+Z/3Zseq.eq} the
second is split by $[1,1,1]-2(1,1,1)$. Hence
\[
\tH_i(\expstwo{5}) =\begin{cases}
                  \integer                         &   i= 10, \\
                  \integer\oplus\integer/4\integer &   i=8, \\
                  \integer/2\integer               &   i=7, \\
                  \integer/3\integer\oplus \integer/3\integer &   i=6, \\
                  0                                & \mbox{else}.
                  \end{cases}
\]

\subsection{The general case}
\label{spherehomology.sec}

We now turn to the general case. As we have seen, the integral homology of
\expstwo{k}\ becomes increasingly complicated, with more and more torsion
arising as $k$ increases. We therefore give a full answer only with
respect to rational co-efficients, where the situation is much 
cleaner, and limit ourselves to calculating only the top three 
integral groups. Not only are these groups computationally tractable
but they also exhibit great regularity.

We begin by calculating the homology of the clipped cube complex.
We will be working mainly with the top end,
where the partitions consist mostly of ones, and for compactness of 
notation we write $\one(k)$ for the partition $(1,\ldots,1)$ of $k$
as a sum of $k$ ones, or just $\one$ if $k$ is understood.
We write $\two{i}(k)$ for the partition 
$(1,\ldots,2,\ldots,1)$ with a single $2$ occurring in the $i$th place,
and write $\three{i}$, $\four{i}$ for the analogous partition with
a single $3$ or $4$. By extension, we
write $\two{i}\three{j}$ for the partition with a $2$ in the
$i$th place, a $3$ in the $j$th place, and the remaining entries
ones.
With these conventions, the top four chain groups of the clipped
cube complex may be written
\begin{align*}
\S_{k,k}   &= \Span\{\one\}  ,              \\
\S_{k-1,k} &= \Span\{\two{i} | 1\leq i\leq k-1\}, \\
\S_{k-2,k} &= \Span\{\three{i}|1\leq i\leq k-2\}\oplus
              \Span\{\two{i}\two{j}| 1\leq i<j\leq k-2\}, \\
\S_{k-3,k} &= \Span\{\four{i}|1\leq i\leq k-3\}\oplus
              \Span\{\two{i}\three{j}|i\not=j, 1\leq i,j\leq k-3\} \\
       &\qquad\qquad\oplus\Span\{\two{i}\two{j}\two{m}|1\leq i<j<m\leq k-3\}.
\end{align*}
Fortunately, there will be no need to count beyond four or
work with juxtapositions of more than two characters.

Since the boundary of $\one$ is zero we see 
immediately that the $k$th homology 
of the clipped cube complex is \integer. The group $H_{k-1}$ is
calculated almost as easily: the boundary of
$\two{i}$ is given by
\[
\partial \two{i} = 
  \begin{cases}
  \three{1}                     & i=1, \\
  (-1)^i(\three{i-1}-\three{i}) & 2\leq i\leq k-2, \\
  \three{k-2}                   & i=k-1, 
  \end{cases}
\]
and taking alternating sums gives
\begin{equation}
\partial \sum_{i=1}^j (-1)^{i-1}\two{i} =
  \begin{cases}
  \three{j} & 1\leq j\leq k-2, \\
  0         & j=k-1. 
  \end{cases}
\label{imageis3.eq}
\end{equation}
This shows that $H_{k-1}$ is isomorphic to $\integer$ also, with
generator ${\displaystyle\sum_{i=1}^{k-1}(-1)^{i-1}\two{i}}$. The
remaining homology groups are rather more difficult to calculate and are
the subject of the following lemma:

\begin{lemma}
The clipped cube complex $(\S_{*,k},\dq)$ is exact over \rational\
at each $\S_{\ell,k}$, $1\leq \ell\leq k-2$, and exact over 
$\integer$ at $\S_{k-2,k}$. 
\label{clippedhomology.lem}
\end{lemma}

\begin{proof}
The proof is in two parts: we first prove exactness over \integer\ at
$\S_{k-2,k}$ directly, and then prove exactness at the remaining 
groups by induction over $k$, using the result for the clipped
cube complexes $(\S_{*,k-1},\dq)$ and $(\S_{*,k-2},\dq)$. The 
direct proof of exactness at $\S_{k-2,k}$ is also required to complete
the induction step. For simplicity, in what follows we write
$\partial$ for $\dq$. 

By~\eqref{imageis3.eq} 
the image of $\partial_{k-1}$ is the
span of $\{\three{i}|1\leq i\leq k-2\}$, so to prove exactness
at $\S_{k-2,k}$ over \integer\ it suffices to show that $\partial_{k-2}$
restricted to the span of $\{\two{i}\two{j}| 1\leq i<j\leq k-2\}$
is injective. Let $U$ be this span and write $\S_{k-3,k}=V\oplus W$, where
\begin{align*}
V &= \Span\{\four{i}|1\leq i\leq k-3\}\oplus
             \Span\{\three{i}\two{j}|1\leq i<j\leq k-3\}, \\
W &= \Span\{\two{i}\three{j}|1\leq i<j\leq k-3\}
       \oplus\Span\{\two{i}\two{j}\two{m}|1\leq i<j<m\leq k-3\}.
\end{align*}
Let $\bar{V}=(V\oplus W)/W$, and observe that
$\rank U =\rank \bar{V} = \binom{k-2}{2}$. 
The boundary induces a map $\db\co U\rightarrow \bar{V}$, 
and we claim that with respect to suitable orderings of the bases
$\{\two{i}\two{j}|i<j\}$ and $\{\four{i}\}\cup\{\three{i}\two{j}|i<j\}$
the matrix of $\db$ is upper triangular.

To see this, lexicographically order each basis according to 
$(j-i,i)$, where for this purpose we regard $\four{i}$ as $\three{i}\two{i}$
and place it in the order as $(0,i)$. The boundary of 
$\two{i}\two{j}$ is given by
\[
\partial\two{i}\two{j} =
 \begin{cases}
 (-1)^{i}(\three{i-1}\two{i}-\four{i}+\two{i}\three{i+1}) 
                    & j=i+1, \\
 (-1)^{i}(\three{i-1}\two{j-1}-\three{i}\two{j-1})
        +(-1)^{j}(\two{i}\three{j-1}-\two{i}\three{j})
                    & j-1\geq 2,
 \end{cases}
\]
in which the terms containing $\three{i-1}$ should be omitted if
$i=1$, and the terms containing $\three{j}$ should be omitted if $j=k-2$. 
Passing to the quotient $\bar{V}$ we have
\[
\db\two{i}\two{j} =
 \begin{cases}
 (-1)^{i}(\three{i-1}\two{i}-\four{i})
                    & j=i+1, \\
 (-1)^{i}(\three{i-1}\two{j-1}-\three{i}\two{j-1})
                    & j-1\geq 2,
 \end{cases}
\]
and in each case we see that the least element occurring in 
$\db\two{i}\two{j}$ is the one ordered by $(j-i-1,i)$.
This implies that the matrix
of $\db$ is lower triangular with each diagonal entry equal to $\pm 1$, 
and we conclude that $\db$ is an isomorphism.

We now show by induction on $k$ that the boundary map
has rank $\binom{k-2}{\ell-2}$ at $\S_{\ell,k}$ for each $1\leq \ell\leq k-2$. 
This is easily verified for the cases $k=3$ and $k=4$ considered in
section~\ref{smallk.sec} and we use these as the base for the
induction. The case $\ell=k-2$ is the one just proved and is required
for the inductive step. For $1\leq \ell\leq k-3$ we relate
$\S_{\ell,k}$ to $\S_{*,k-1}$ and $\S_{*,k-2}$ 
by examining the size of the last element of each partition. 
Let
\[
\S^i_{\ell,k} = \Span\{S\in\S_{\ell,k}|s_\ell =i\}
\]
for $i=1,2$,
\[
\S^3_{\ell,k} = \Span\{S\in\S_{\ell,k}|s_\ell \geq 3\},
\]
and consider the ``append'' and ``plus'' operators
\begin{align*}
A_i(S) &= (s_1,\ldots,s_\ell,i) \\
\intertext{for $i=1,2$, and}
P_2(S) &= (s_1,\ldots,s_\ell +2).
\end{align*}
We have
\begin{align*}
A_i\co&\S_{\ell-1,k-i}\decarrow{\cong}\S^i_{\ell,k}, \\
P_2\co&\S_{\ell,k-2}\decarrow{\cong}\S^3_{\ell,k},
\end{align*}
and moreover
\[
\partial A_i S \in A_i\partial S + \S^3_{\ell,k}
\]
for $i=1,2$, and
\begin{equation}
\label{plus2commutes.eq}
\partial P_2 S =   P_2\partial S
\end{equation}
since $P_2$ does not change the parity of the last entry of $S$.
We note also that~\eqref{plus2commutes.eq} depends on the fact we are
using the boundary map $\dq$ instead of $\dlambda$.

Since $\S_{\ell,k}=\S^1_{\ell,k}\oplus\S^2_{\ell,k}\oplus\S^3_{\ell,k}$ it
follows that the matrix of $\partial_\ell$ can be written in the block
form
\[
D_{\ell,k} =
\begin{bmatrix}
D_{\ell-1,k-1} &  0              &    0     \\
0              &  D_{\ell-1,k-2} &    0  \\
E              &  F              &    D_{\ell,k-2} 
\end{bmatrix},
\]
in which $D_{i,j}$ is the matrix of $\partial$ acting on $\S_{i,j}$. 
Hence
\begin{align*}
\rank \partial_\ell &= 
          \rank D_{\ell-1,k-1} + \rank D_{\ell-1,k-2} + \rank D_{\ell,k-2} \\
       &= \binom{k-3}{\ell-3}  + \binom{k-4}{\ell-3}  +\binom{k-4}{\ell-2} \\
       &= \binom{k-2}{\ell-2}
\end{align*}
as desired.

We complete the proof using a rank argument. We have 
\begin{align*}
\rank \ker \partial_\ell &= \rank \S_{\ell,k} - \rank\partial_\ell \\
                         &= \binom{k-1}{\ell-1} - \binom{k-2}{\ell-2} \\
                         &= \binom{k-2}{\ell-1} \\
                         &  =\rank \partial_{\ell+1},
\end{align*}
from which the result follows.
\end{proof}

Since $\dq$ and $\dlambda$ co-incide on $\S_{k,k}$ and $\S_{k-1,k}$ we
may restate this lemma and the calculations preceding it in the 
following form:

\begin{lemma}
The rational homology of the space $\expstwob{k}/\expstwob{k-1}$
vanishes except in dimensions $2k-2$ and $2k-3$. With respect to
integer co-efficients the top three groups are \integer\ in
dimension $2k-2$, generated by $\one(k-1)$, \integer\ in dimension $2k-3$,
generated by ${\displaystyle\sum_{i=1}^{k-2}(-1)^{i-1}\two{i}(k-1)}$, 
and $0$ in dimension $2k-4$.
\end{lemma}

We are now in a position to prove Theorems~\ref{s2rational.th} 
and~\ref{s2integral.th}. The argument should be very familiar 
from the cases considered in section~\ref{smallk.sec}.

\begin{proof}[Proof of Theorems~\ref{s2rational.th} and~\ref{s2integral.th}]
We first prove both results for $\expstwob{k}$, using induction on $k$.
The top end of the long exact sequence of the pair
$(\expstwob{k},\expstwob{k-1})$ gives 
\[
0\longrightarrow  \tH_{2k-2}(\expstwob{k}) \longrightarrow 
       \tH_{2k-2}(\expstwob{k}/\expstwob{k-1})  \longrightarrow 0,
\]
implying $\tH_{2k-2}(\expstwob{k})\cong\integer$ with generator
$\one(k-1)$. The next segment is 
\[
0\longrightarrow  \tH_{2k-3}(\expstwob{k}) \longrightarrow 
       \integer \decarrow{\partial}
       \integer  \longrightarrow 
       \tH_{2k-4}(\expstwob{k})\longrightarrow 0,
\]
where the boundary map in the middle, from
$\tH_{2k-3}(\expstwob{k}/\expstwob{k-1})$ to 
$\tH_{2k-4}(\expstwob{k-1})$, is given by $\dnu$. Since
$\dnu\two{i}(k-1)=(-1)^{i}\one(k-2)$ this sends the
generator ${\displaystyle\sum_{i=1}^{k-2}(-1)^{i-1}\two{i}(k-1)}$
to $(2-k)(\one(k-2))$, inducing multiplication by $2-k$. It follows
that $\tH_{2k-3}(\expstwob{k})$ is zero 
and that $\tH_{2k-4}(\expstwob{k})$ is isomorphic to 
$\integer/(k-2)\integer$, generated by the top homology class of
\expstwob{k-1}. We show that the remaining groups vanish over
\rational\ using the induction hypothesis and the exact sequence
\[
\tH_i(\expstwob{k-1})\longrightarrow
\tH_i(\expstwob{k})\longrightarrow
\tH_i(\expstwob{k}/\expstwob{k-1})
\]
(rational co-efficients omitted). 
Since the outer groups are zero for $i\leq 2k-5$ the middle group is too.

We now prove the results for \expstwo{k}, using the results just
proved, Lemma~\ref{samehomology.lem}, and the long exact sequence
of the pair $(\expstwo{k},\expstwob{k})$. At the top end
we have
\[
0\rightarrow  \tH_{2k}(\expstwo{k}) \rightarrow 
       \tH_{2k}(\expstwo{k}/\expstwob{k})  \rightarrow
0\rightarrow \tH_{2k-1}(\expstwo{k})\rightarrow 0
\]
giving $\tH_{2k}(\expstwo{k})\cong\integer$ generated by
$\tilde\one(k)$ and $\tH_{2k-1}(\expstwo{k})\cong\{0\}$. The next
segment is 
\[
0\longrightarrow \integer \longrightarrow \tH_{2k-2}(\expstwo{k})
 \longrightarrow \integer/(k-1)\integer \longrightarrow 0
\]
in which the subgroup \integer\ is generated by $\one(k-1)$ and
the quotient $\integer/(k-1)\integer$ is generated by $\tilde\one(k-1)$. 
This short exact sequence splits: the cycle $\tilde\one(k-1)-2(\one(k-1))$ maps
to $\tilde\one(k-1)$ in the quotient, and
\[
\partial \sum_{i=1}^{k-1}(-1)^{i-1}\ttwo{i}(k)
           =(k-1)\bigl(\tilde\one(k-1)-2(\one(k-1))\bigr),
\]
showing that this class has order $k-1$. This completes the proof of
Theorem~\ref{s2integral.th}, and the exact sequence
\[
\tH_i(\expstwob{k};\rational)\longrightarrow
\tH_i(\expstwo{k};\rational)\longrightarrow
\tH_i(\expstwo{k}/\expstwob{k};\rational)
\]
yields $0\longrightarrow \tH_i(\expstwo{k};\rational)\longrightarrow 0$ for
$i\leq 2k-3$, completing the proof of Theorem~\ref{s2rational.th}.
\end{proof}

\section{Higher genus surfaces}
\label{highergenus.sec}

In this section we prove Theorem~\ref{surfacehomology.th}
on the top homology of the finite subset spaces of a closed surface,
and Theorem~\ref{surfacedegree.th} on the degree of a
map $\Exp{k}{f}$
induced
by a map $f\co\Sigma\rightarrow\Sigma'$ between closed oriented
surfaces. 

We give $\Sigma$ the standard cell structure consisting of
a single vertex, $2g$ or $g$ edges, and a single two-cell 
attached along the word $w_+=[e_1,e_{1+g}]\cdots[e_g,e_{2g}]$
if $\Sigma$ is orientable, $w_-=e_1^2\cdots e_g^2$ if $\Sigma$ is
non-orientable. With respect to this cell structure the $2k$-- and
$(2k-1)$--cells of the lexicographic
cell structure for $\expsig{k}$ are the single cell $\tilde\one(k)$
in dimension $2k$ and the cells
\[
\bigl\{\ttwo{i}(k)|1\leq i\leq k-1\bigr\}\cup
\bigl\{(\tilde\one(k-1))\cup\ecell[i]{1}|
                1\leq i\leq \dim H_1(\Sigma;\integer/2\integer)\bigr\}
\]
in dimension $2k-1$. In order to calculate boundaries we need
to know the chain maps $(\expinf{w_+})_\sharp$ and 
$(\expinf{w_-})_\sharp$. By the discussion at the end of
section~3.3 of~\cite{graphs03} it suffices to calculate the
images of the generators $\ecell{1}$ and $\ecell{2}$ of
the chain ring, and moreover these are
given by
\begin{align*}
(\expinf{w_+})_\sharp\ecell{1} &= 0,                     &
(\expinf{w_-})_\sharp\ecell{1} &= 2\sum_{i=1}^g \ecell[i]{1}, \\
(\expinf{w_+})_\sharp\ecell{2} &= 
           2\sum_{i=1}^g \ecell[i]{1}\cup\ecell[i+g]{1}, &
(\expinf{w_-})_\sharp\ecell{2} &= 2\sum_{i=1}^g \ecell[i]{2}
           +4\sum_{i<j} \ecell[i]{1}\cup\ecell[j]{1}.
\end{align*}
This gives us all the ingredients we need and we now proceed to the proofs.
\begin{proof}[Proof of Theorem~\ref{surfacehomology.th}]
Since $\dnu\tilde\one=\dlambda\tilde\one=0$ the boundary of the cell 
$\tilde\one(k)$ is 
equal to $\pm\dgamma\tilde\one(k)$. We have
\begin{equation}
\dgamma\tilde\one(k) = 
 \begin{cases}
 0                                       & \text{$\Sigma$ orientable}, \\
 {\displaystyle 2\sum_{i=1}^g \tilde\one(k-1)\cup\ecell[i]{1}}
                                         & \text{$\Sigma$ non-orientable},
 \end{cases}
\label{dgammaone.eq}
\end{equation}
and we see immediately that
\[
H_{2k}(\expsig{k})=\begin{cases}
                    \integer & \mbox{if $\Sigma$ is orientable,} \\
                    0    & \mbox{if $\Sigma$ is non-orientable.}
                   \end{cases}
\]

To calculate the kernel of $\partial_{2k-1}$ we make the following
observation. Let $F_d$ be the span of the $d$--cells of the form
$\vfcell{S}$, $\fcell{S}$, and let $E_d$ be the span of the 
$d$--cells with a $\vcell{J}$ or $\ecell{J}$ factor with $|J|>0$.
With respect to the decompositions of the cellular chain groups
$\mathcal{C}_d$ as $F_d\oplus E_d$ the
matrix of $\partial$ has the block form
\[
D_d = \twobytwo{A_d}{0}{B_d}{C_d},
\]
and if $A_d$ is injective then $\ker\partial =\ker C_d$. Applying this
with $d=2k-1$ we have
\begin{align*}
F_{2k-1} &= \Span\{\ttwo{i}(k)|1\leq i\leq k-1\}, \\
E_{2k-1} &= \Span\{\tilde\one(k-1)\cup
      \ecell[i]{1}|1\leq i\leq \dim H_1(\Sigma;\integer/2\integer)\},
\end{align*}
and $A_{2k-1}$ is injective since 
$A_{2k-1}\ttwo{i}$ is equal to the boundary of this cell as a cell
in \expstwo{k}. It follows that 
the $(2k-1)$--cycles are contained in the span of 
$\{\tilde\one(k-1)\cup\ecell[i]{1}|1\leq i\leq 
  \dim H_1(\Sigma;\integer/2\integer)\}$.
By equation~\eqref{dgammaone.eq} we have
\[
\partial\bigl(\tilde\one(k-1)\cup\ecell[i]{1}\bigr) = 
 \begin{cases}
 0                                       & \text{$\Sigma$ orientable}, \\
 {\displaystyle\pm2\sum_{j\not=i} 
     \tilde\one(k-2)\cup\ecell[j]{1}\cup\ecell[i]{1} }
                                         & \text{$\Sigma$ non-orientable}.
 \end{cases}
\]
If $\Sigma$ is non-orientable then 
\[
\partial\sum_{i=1}^g a_i\tilde\one(k-1)\cup\ecell[i]{1} =
   \pm 2\sum_{i<j} (a_i-a_j)\tilde\one(k-2)\cup\ecell[i]{1}\cup\ecell[j]{1},
\]
implying
\[
\ker\partial_{2k-1} =
 \begin{cases}
\Span\bigl\{\tilde\one(k-1)\cup\ecell[i]{1}|1\leq i\leq 2g\bigr\}
                                       & \text{$\Sigma$ orientable},\\
{\displaystyle
      \Span\biggl\{\sum_{i=1}^g \tilde\one(k-1)\cup\ecell[i]{1}\biggr\}}
                                       & \text{$\Sigma$ non-orientable}.
 \end{cases}
\]
It now follows that
\[
H_{2k-1}(\expsig{k})=\begin{cases}
                    \integer^{2g} & \mbox{if $\Sigma$ is orientable,} \\
                    \integer/2\integer & \mbox{if $\Sigma$ is non-orientable,}
                    \end{cases}
\]
completing the proof of Theorem~\ref{surfacehomology.th}.
\end{proof}

Next we prove Theorem~\ref{surfacedegree.th}, which we recall states that
\[
\deg \Exp{k}{f} = (\deg f)^k
\] 
if $f\co\Sigma\rightarrow\Sigma'$ is a map between closed oriented
surfaces. The result is
an almost immediate consequence of the fact that 
$\expsig{k}\setminus\expsig{k-1}$ may be consistently oriented
by~\eqref{surfaceorientation.eq}. We caution that this orientation
disagrees with the orientation on $\tilde\one(k)$ for some $k$.

\begin{proof}[Proof of Theorem~\ref{surfacedegree.th}]
Given a map
$f\co\Sigma\rightarrow\Sigma'$ let 
$p_1,\ldots,p_k$ be distinct points in $\Sigma'$ and perturb
$f$ to be transverse to each $p_i$. Each point $p_i$ then has
finitely many preimages which we label
$q_{i,1},\ldots,q_{i,r_i}$.
Let $\Lambda=\{p_1,\ldots,p_k\}\in \Exp{k}{\Sigma'}$. Under
\Exp{k}{f}\ the preimage of $\Lambda$ consists of the
$r_1\cdots r_k$ points
\[
\Lambda_{s_1\cdots s_k} = \{q_{1,s_i},\ldots,q_{k,s_k}\},
\]
and perturbing the characteristic maps of the $2$--cells of 
$\Sigma$, $\Sigma'$ if necessary we may assume that all points in
question lie in the top dimensional cells $\tilde\one(k)$. With
respect to the canonical orientation on $\expsig{k}\setminus\expsig{k-1}$ 
the sign
of $\Lambda_{s_1\cdots s_k}$ is the product $\prod_i\sign q_{i,s_i}$, 
and we have immediately
\[
\deg\Exp{k}{f} = \prod_{i=1}^k \sum_{j=1}^{r_i} \sign q_{i,j} = (\deg f)^k.
\]
\end{proof}

\section{A homotopy model for the third finite subset space}
\label{model.sec}

We now pursue our second direction in the study of the finite subset
spaces of closed surfaces. We begin by constructing a homotopy model
for the third finite subset space of a closed orientable surface, and
then we use this in section~\ref{exp3sigcohomology.sec} to calculate
the cohomology of $\expsig{3}$. 

\subsection{The model}

To construct a homotopy model for \expsig{3}\ we begin with the 
intermediate quotient $\symthree =\Sigma^3/S_3$. Inside \symthree\ we
have the preimage of \expsig{2}, namely the branch locus \branch\
consisting of the quotient of the diagonals of 
$\Sigma\times\Sigma\times\Sigma$. Instead of quotienting \branch\ 
further to get \expsig{3}\ we attach the mapping cylinder of 
$\branch\rightarrow\expsig{2}$
to obtain a homotopy equivalent space \model . Fortunately \branch\ is 
simply a copy of $\Sigma\times\Sigma$ and $\branch\rightarrow\expsig{2}$
is the quotient map $\Sigma\times\Sigma\rightarrow\symtwo$, so the 
construction is all in terms of known spaces and maps.

Concretely, let $q_k\co\Sigma^k\rightarrow\sym{k}{\Sigma}$ be the quotient
map and $\Delta\co\Sigma\rightarrow\Sigma\times\Sigma$ the diagonal map. 
Then \branch\ is the image of $\Sigma\times\Sigma$ under the
map $\iota=q_3 \circ\id\times\Delta$, and is homeomorphic to 
$\Sigma\times\Sigma$ since $\iota$ is injective, $\Sigma\times\Sigma$ is
compact, and \symthree\ is Hausdorff. Further, \expsig{3}\ is
obtained from \symthree\ by identifying $\iota(a,b)=q_3(a,b,b)$ with
$\iota(b,a)=q_3(b,a,a)$. Let 
\[
M_{q_2} = \frac{(\Sigma\times\Sigma\times I) \amalg \symtwo}
               {(x,1)\sim q_2(x)}
\]
be the mapping cylinder of $q_2$, and note that we adopt the convention that 
mapping cylinders are attached to the target at $1\in I$ rather than $0$.
We obtain our model for \expsig{3}\ by attaching $M_{q_2}$ to 
\symthree\ along \branch\ and $\Sigma\times\Sigma\times\{0\}$, namely
\begin{align*}
\model &=  \symthree \cup_{\Sigma\times\Sigma} M_{q_2} \\
       &=  \frac{\symthree \amalg M_{q_2}}{\iota(x)\sim (x,0)} .
\end{align*}

\begin{figure}[t]
\begin{center}
\leavevmode
\psfrag{sym3}{\symthree}
\psfrag{sym2}{\symtwo}
\psfrag{D=Sx0}{$\branch=\Sigma^2\!\times\!\{0\}$}
\psfrag{Sx1/2}{$\Sigma^2\!\times\!\{1/2\}$}
\psfrag{Mq}{$M_{q_2}$}
\includegraphics{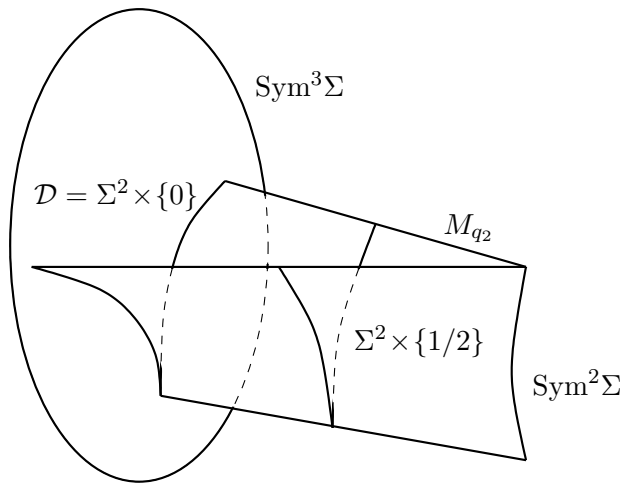}
\caption[A homotopy model for the third finite subset space
of a surface]
{A schematic picture of the space $\model$. As we shall see
in section~\protect\ref{SymkS.sec}, \symthree\ 
is a complex $3$--manifold with the branch locus $\branch\cong\Sigma^2$
embedded with a cusp along the diagonal. $\Sigma^2$ maps two-to-one to 
\symtwo\ via $q_2$, and we form $\model$ by attaching the mapping cylinder 
$M_{q_2}$ to $\branch$ along $\Sigma^2\times\{0\} $.}
\label{model.fig}
\end{center}
\end{figure}

The space \model\ is shown schematically in figure~\ref{model.fig}. The 
subscript $3$ reflects the hope that a similar construction may apply for
$k\geq 4$, perhaps using several mapping cylinders to successively quotient
\sym{k}{\Sigma}\ to \expsig{k}\ in several stages. However, such a 
generalisation is complicated by the increasing complexity of the
branch locus: for example, when $k=4$ it is 
$\Sigma\times\symtwo$ with an embedded copy of $\Sigma\times\Sigma$
quotiented to $\symtwo$.

The cornerstone of the next section is the following lemma:

\begin{lemma}
The spaces \model\ and \expsig{3}\ are homotopy equivalent.
\label{homotopic.lem} 
\end{lemma}

\noindent
We prove the lemma after a brief digression on known facts about
symmetric products of
surfaces. Not all of what follows is central to our argument, but 
it nonetheless serves to give a fuller picture of the construction.

\subsection{Symmetric products of surfaces}
\label{SymkS.sec}

We recall that the symmetric product of an orientable surface is
a manifold, and moreover that a complex structure on $\Sigma$ leads
to a complex structure on $\sym{k}{\Sigma}$.  Local co-ordinates about 
$\{p_1,\ldots,p_k\}$ are given by the elementary symmetric functions 
in the local complex co-ordinates about the $p_i$, and are obtained
by regarding the points as the zeroes of a polynomial (see
Griffiths and Harris~\cite[p.\ 326]{griffiths-harris-poag}). In particular
$\sym{n}{\complex P^1}$ may be identified with the non-vanishing homogeneous
polynomials of degree $n$ in two variables, modulo scaling, and as such
is equal to $\complex P^n$. 

A polynomial $p$ has repeated roots if and only if its discriminant
is zero. The discriminant is a polynomial in the co-efficients of $p$
(see Lang~\cite[pp.\ 192--194]{lang-algebra})
and it follows that the branch locus, the image of the diagonals of
$\Sigma^k$, is locally given by the vanishing of a polynomial and 
is therefore an algebraic variety. Specialising to $k=3$, for a suitable choice
of complex co-ordinates $(a,b,c)$ the branch locus \branch\ is locally the set
$b^2=c^3$. We see that \branch\ has a cusp along $b=c=0$, which is precisely
the image of the main diagonal in $\Sigma\times\Sigma\times\Sigma$.

\subsection{The proof of homotopy equivalence}

Consider the mapping cylinder of $q\co\symthree\rightarrow\expsig{3}$. This
deformation retracts to \expsig{3}\ and contains \model\ as a subspace, 
and our aim is to show that it also deformation retracts to \model. 
To do this it suffices to show that $\symthree\times I$ deformation 
retracts to $\symthree\times\{0\}\cup\branch\times I$, as such a 
homotopy will descend setwise to the quotient $M_q$ and be continuous
there.

The existence of a deformation retraction from $\symthree\times I$ to
$\symthree\times\{0\}\cup\branch\times I$
follows from results in Bredon~\cite[pp.\ 431--432]{bredon} and
Dugundji~\cite[pp.\ 327--328]{dugundji} and the existence of a 
neighbourhood $U\supseteq\branch$ that strongly deforms to \branch\ in
\symthree. We will prove it in this
way, using the fact that \symthree\ and \branch\ are both compact
manifolds to construct $U$. However, an approach of
perhaps greater generality might be to realise \branch\ as a subcomplex
of \symthree\ and appeal to Hatcher~\cite[Prop.\ 0.16]{hatcher}.
Hatcher~\cite[pp.\ 482--483]{hatcher} gives a construction of an
$S_k$--equivariant simplicial structure on $X^k$ for simplicial $X$, 
and Lemma~\ref{homotopic.lem} would follow from checking whether this
contains the diagonals as a subcomplex. 

Although \symthree\ and \branch\ are both manifolds the existence of 
the desired neighbourhood $U$ does not simply follow
from the tubular neighbourhood theorem, since \branch\ is not smoothly
embedded. We therefore resort to more hands-on means, and use the fact
that manifolds are Euclidean neighbourhood retracts. Embed $M=\symthree$
in some $\real^n$. Then there are neighbourhoods $V$ of $M$, $W$ of \branch, 
and retractions $r_M\co V\rightarrow M$, $r_W\co W\rightarrow \branch$. 
The neighbourhood
$W$ may be taken sufficiently small that the linear homotopy from 
$W\hookrightarrow\real^n$ to $r_W$ remains in $V$, and we post-compose 
this with $r_M$ and intersect $W$ with $M$ to get the desired neighbourhood 
and deformation.

\section{The calculation of $H^*(\expsig{3})$}
\label{exp3sigcohomology.sec}

\subsection{Introduction}

To calculate the cohomology of \model\ we use the Mayer-Vietoris
sequence and the obvious decomposition
\[
\model = (\model\setminus\symtwo)\cup(\model\setminus\symthree).
\]
The pieces are homotopy equivalent to \symthree\ and \symtwo\
respectively, and intersect in 
$\Sigma\times\Sigma\times(0,1)\simeq\Sigma\times\Sigma$, leading to  
a long exact sequence
\[
\cdots\rightarrow H^i(\model)\rightarrow H^i(\symthree)\oplus H^i(\symtwo)
\rightarrow H^i(\Sigma^2)\rightarrow H^{i+1}(\model)\rightarrow\cdots .
\]
Before proceeding we describe the rings $H^*(\Sigma^k)$ and 
$H^*(\sym{k}{\Sigma})$. Recall that integer co-efficients are to be
assumed except where specified otherwise.

\subsection{The rings $H^*(\Sigma^k)$ and $H^*(\sym{k}{\Sigma})$}

Let $\alpha_1,\ldots,\alpha_{2g}$ be generators for $H^1(\Sigma)$ such 
that 
\[
\alpha_i\alpha_j =
\begin{cases}
0 & |i-j|\not=g, \\
\beta & j=i+g ,
\end{cases}
\]
where $\beta$ is a generator of $H^2(\Sigma)$. Since $H^*(\Sigma)$ is
finitely generated and free the K\"unneth formula applies and
\[
H^*(\Sigma^k) \cong H^*(\Sigma)^{\otimes k}.
\]

The cohomology ring of \sym{k}{\Sigma}\ is given by
Macdonald~\cite{macdonald62} and 
Seroul~\cite{seroul72-pdm}. In addition Seroul's paper~\cite{seroul72-crasci}
gives a sketch of his argument. Macdonald
uses methods from algebraic geometry to give generators and relations
for $H^*(\sym{k}{\Sigma};K)$ over a field $K$ of characteristic zero,
and to show that $H^*(\sym{k}{\Sigma})$ is torsion free. He then 
states incorrectly that this implies the same elements generate over
the integers. Seroul confirms Macdonald's answer, using purely 
algebraic-topological techniques to find
$H^*(\sym{k}{\Sigma};\integer)$ directly.
In part the result may be stated as follows; we omit the statement
of the relations as we will do all ring multiplication in 
$H^*(\sym{k}{\Sigma})$. 

\begin{theorem}[Macdonald~\cite{macdonald62} and 
Seroul~\cite{seroul72-pdm,seroul72-crasci}]
The map
\[
q_k^*\co H^*(\sym{k}{\Sigma};R)\rightarrow H^*(\Sigma^k;R)
\]
is an isomorphism of $H^*(\sym{k}{\Sigma};R)$ onto $H^*(\Sigma^k;R)^{S_k}$,
the subring of cohomology fixed by $S_k$, for $R$ a field of characteristic
zero, and is injective for $R=\integer$. $H^*(\sym{k}{\Sigma};\integer)$
is generated by elements $\xi_1,\ldots,\xi_{2g}$ in degree $1$ and
$\eta$ in degree $2$ such that
\begin{align*}
q_k^*\xi_i &= \sum_{j=1}^k \pi_j^*\alpha_i, & 
q_k^*\eta &= \sum_{j=1}^k \pi_j^*\beta,
\end{align*}
where $\pi_j\co\Sigma^k\rightarrow\Sigma$ is projection on the $j$th factor.
A basis for $H^r(\sym{k}{\Sigma};\integer)$ is given by the monomials
$\xi_{i_1}\cdots\xi_{i_m}\eta^n$  for which $m+2n=r$, 
$i_1<\cdots<i_m$, and $m\leq\min\{r,2k-r\}$.
\end{theorem}
\noindent
We remark that $q_k$ is a degree $k!$ map, so $q_k^*$ is certainly
not onto $H^*(\Sigma^k)^{S_k}$ with integer co-efficients.

To avoid confusion we give different names to the generators of 
$H^*(\symtwo)$ and $H^*(\symthree)$. Let
\begin{gather*}
\zeta_i = \alpha_i\otimes 1 + 1\otimes\alpha_i, \\
\theta = \beta\otimes 1 + 1\otimes\beta ,        \\
\xi_i = \alpha_i\otimes 1\otimes1 + 1\otimes\alpha_i\otimes 1 
            + 1\otimes 1\otimes\alpha_i ,     \\
\eta = \beta\otimes 1\otimes1 + 1\otimes\beta\otimes 1 
            + 1\otimes 1\otimes\beta . 
\end{gather*}
Since $q_k^*$ is injective we shall abuse notation and not take
care to distinguish between elements of $H^*(\sym{k}{\Sigma})$ and
their images in $H^*(\Sigma^k)$, and will regard
the  $\zeta_i$ and $\theta$ as generators of $H^*(\symtwo)$, and 
the $\xi_i$ and $\eta$ as generators of $H^*(\symthree)$. 

\subsection{The cohomology calculation}

Returning to the Mayer-Vietoris sequence, letting
\[
\Phi_i\co H^i(\symthree)\oplus H^i(\symtwo)\rightarrow H^i(\Sigma\times\Sigma)
\]
be the map $\iota^*\oplus q_2^*=(q_3\circ\id\times\Delta)^*\oplus q_2^*$
we have the short exact sequence
\[
0\rightarrow \coker \Phi_{i-1} \rightarrow H^i(\model)\rightarrow
\ker\Phi_i\rightarrow 0.
\]
Since $H^i(\symthree)\oplus H^i(\symtwo)$ is free the kernel of $\Phi_i$ 
is too, so the sequence splits and we get
\begin{equation}
\label{directsum.eq}
H^i(\model)\cong\coker\Phi_{i-1}\oplus\ker\Phi_i .
\end{equation}
In what follows we calculate the kernel and cokernel of each $\Phi_i$.

\subsubsection*{Dimension one}

Both $H^1(\symthree)\oplus H^1(\symtwo)$ and $H^1(\Sigma\times\Sigma)$ have
rank $4g$, with bases $\{\xi_i\}\cup\{\zeta_i\}$ and 
$\{\alpha_i\otimes 1\}\cup\{1\otimes\alpha_i\}$ respectively. For 
$\alpha\in H^j(\Sigma)$ we have
\begin{align*}
(\id\times\Delta)^*(1\otimes 1\otimes\alpha)
                    &= (\id\times\Delta)^*\pi_3^*\alpha \\
                    &= (\pi_3\circ \id\times\Delta)^*\alpha \\
                    &= \pi_2^*\alpha = 1\otimes\alpha,
\end{align*}
and similarly $(\id\times\Delta)^*(1\otimes\alpha\otimes 1)=1\otimes\alpha$,
$(\id\times\Delta)^*(\alpha\otimes 1\otimes 1)= \alpha\otimes 1$.
Consequently
\begin{align*}
\Phi_1(\xi_i) &= (\id\times\Delta)^*(\alpha_i\otimes 1\otimes1 
                 + 1\otimes\alpha_i\otimes 1 + 1\otimes 1\otimes\alpha_i) \\
                &= \alpha_i\otimes 1 + 2\otimes\alpha_i .
\end{align*}
Since $\Phi_1(\zeta_i)=q_2^*\zeta_i = \alpha_i\otimes 1 + 1\otimes\alpha_i$
and $\det\twobytwo{1}{1}{2}{1}=-1$, $\Phi_1$ maps $\Span\{\xi_i,\zeta_i\}$
isomorphically onto $\Span\{\alpha_i\otimes 1,1\otimes\alpha_i\}$. Thus
\[
\ker\Phi_1 \cong\coker\Phi_1\cong\{0\}.
\]
\subsubsection*{Dimension two}

$H^2(\symthree)\oplus H^2(\symtwo)$ has basis
\[
\{\xi_i\xi_j|i<j\}\cup\{\zeta_i\zeta_j|i<j\}\cup\{\eta,\theta\}
\]
and rank $2\binom{2g}{2}+2$, while $H^2(\Sigma\times\Sigma)$ has basis
\[
\{\alpha_i\otimes\alpha_j\}\cup\{\beta\otimes 1,1\otimes\beta\}
\]
and rank $4g^2+2$. Under $\Phi_2$ we have
\begin{align*}
\xi_i\xi_j &\mapsto (\alpha_i\otimes 1 + 2\otimes\alpha_i)
                        (\alpha_j\otimes 1 + 2\otimes\alpha_j) \\
     &= \begin{cases}
     2(\alpha_i\otimes\alpha_j  - \alpha_j\otimes\alpha_i) & |i-j|\not=g, \\
     2(\alpha_i\otimes\alpha_j  - \alpha_j\otimes\alpha_i) 
               +\beta\otimes 1 +4\otimes\beta & j=i+g,
        \end{cases}   \\
\zeta_i\zeta_j &\mapsto (\alpha_i\otimes 1 + 1\otimes\alpha_i)
                    (\alpha_j\otimes 1 + 1\otimes\alpha_j) \\
     &= \begin{cases}
     \alpha_i\otimes\alpha_j  - \alpha_j\otimes\alpha_i & |i-j|\not=g, \\
     \alpha_i\otimes\alpha_j  - \alpha_j\otimes\alpha_i
               +\beta\otimes 1 +1\otimes\beta & j=i+g,
     \end{cases} \\
\eta &\mapsto \beta\otimes 1 + 2\otimes \beta, \\
\theta &\mapsto \beta\otimes 1 + 1\otimes \beta.
\end{align*}
Clearly the image of $\Phi_2$ is the span of
\begin{equation}
\{\beta\otimes 1,1\otimes\beta\}
\cup\{\alpha_i\otimes\alpha_j-\alpha_j\otimes\alpha_i|i<j\},
\label{BasisForImPhi_2.eq}
\end{equation}
a subspace of rank $\binom{2g}{2}+2$. Thus the kernel of $\Phi_2$ has 
rank $\binom{2g}{2}$. The set in~\eqref{BasisForImPhi_2.eq}
may be augmented to a basis for $H^2(\Sigma\times\Sigma)$, so the
cokernel of $\Phi_2$ is free of rank 
$4g^2+2-\binom{2g}{2}-2=\binom{2g}{2}+2g$. Hence
\begin{align*}
\ker\Phi_2 & \cong \integer^{\binom{2g}{2}}, & 
\coker\Phi_2 &\cong \integer^{\binom{2g}{2}+2g}.
\end{align*}

\subsubsection*{Dimension three}

A basis for $H^3(\symthree)\oplus H^3(\symtwo)$ is given by
\[
\{\xi_i\xi_j\xi_k|i<j<k\}\cup\{\xi_i\eta\}\cup\{\zeta_i\theta\}.
\]
If the genus of $\Sigma$ is greater than one the rank 
is $\binom{2g}{3}+4g$, but in genus equal to one there are only 
two distinct $\xi_i$, so the leftmost set in this union is empty 
and the rank of $H^3(\symthree)\oplus H^3(\symtwo)$ is $4g=4$. In 
either case $H^3(\Sigma\times\Sigma)$ has basis
\[
\{\alpha_i\otimes\beta\}\cup\{\beta\otimes\alpha_i\}
\]
and rank $4g$. We have
\begin{align*}
\zeta_i\theta & \mapsto (\alpha_i\otimes 1 + 1\otimes\alpha_i)
                   (\beta\otimes 1 + 1\otimes\beta) \\
         &= \alpha_i\otimes\beta + \beta\otimes\alpha_i, \\
\xi_i\eta &\mapsto (\alpha_i\otimes 1 + 2\otimes\alpha_i)
                   (\beta\otimes 1 +2\otimes\beta) \\
         &= 2(\alpha_i\otimes\beta + \beta\otimes\alpha_i),
\end{align*}
and in genus one it follows that the kernel and cokernel of $\Phi_3$
both have rank two. When $g\geq 2$ the triple product 
$\xi_i\xi_j\xi_k$ maps to $0$ if $i,j,k$ are distinct mod $g$, 
while 
\begin{align*}
\xi_i\xi_{i+g}\xi_j&\mapsto 
  (2(\alpha_i\otimes\alpha_{i+g}  - \alpha_{i+g}\otimes\alpha_i)+
      \beta\otimes 1 +4\otimes\beta)(\alpha_j\otimes 1+2\otimes\alpha_j) \\
  &= 2\beta\otimes\alpha_j+4\alpha_j\otimes\beta
\end{align*}
for $i\not=j\not= i+g$. Considering the images of $\zeta_i\theta$ and 
$\xi_i\xi_{i+g}\xi_j$ we see that the image of $\Phi_3$ has rank
$4g$ and that 
\[
\coker \Phi_3 \cong \frac{\Span\{\beta\otimes\alpha_j+2\alpha_j\otimes\beta\}}
{\Span\{2(\beta\otimes\alpha_j+2\alpha_j\otimes\beta)\}} 
\cong [\integer/2\integer]^{2g},
\]
so that 
\begin{align*}
\ker\Phi_3 &\cong\begin{cases}
                 \integer^2 & g=1, \\
                 \integer^{\binom{2g}{3}} & g\not=1,
                 \end{cases} &
\coker\Phi_3 &\cong \begin{cases}
                    \integer^2 & g=1, \\
                    [\integer/2\integer]^{2g} & g\not=1.
                    \end{cases} 
\end{align*}

\subsubsection*{Dimension four}

$H^4(\symthree)\oplus H^4(\symtwo)$ has rank $\binom{2g}{2}+2$ and basis
\[
\{\xi_i\xi_j\eta|i<j\}\cup\{\eta^2,\theta^2\},
\]
while $H^4(\Sigma\times\Sigma)$ has rank one and basis $\{\beta\otimes\beta\}$.
Under $\Phi_4$ we have
\begin{align*}
\theta^2 &\mapsto (\beta\otimes 1 + 1\otimes\beta)^2 \\
      &= 2\beta\otimes\beta, \\
\eta^2&\mapsto (\beta \otimes 1 + 2\otimes\beta)^2 \\
      &= 4\beta\otimes\beta, \\
\xi_i\xi_j\eta &\mapsto 2(\alpha_i\otimes\beta + \beta\otimes\alpha_i)
                             (\alpha_j\otimes 1 + 2\otimes\alpha_j) \\
      &= \begin{cases}
         0                   & |j-i|\not= g, \\
         6\beta\otimes\beta  & j=i+g.
         \end{cases}
\end{align*} 
Clearly
\begin{align*}
\ker\Phi_4 & \cong \integer^{\binom{2g}{2}+1}, & 
\coker\Phi_4 &\cong \integer/2\integer.
\end{align*}

\subsubsection*{Dimensions five and six}

$\Sigma\times\Sigma$ and $\symtwo$ have no cohomology in dimensions
five and six so the cokernel of $\Phi_i$ is trivial and the kernel
is $H^i(\symthree)$ for $i=5,6$. $H^5(\symthree)=\Span\{\xi_i\eta^2\}$
has rank $2g$ and $H^6(\symthree)=\Span\{\eta^3\}$ has rank one, so
\begin{align*}
\ker\Phi_5 & \cong \integer^{2g}, & \coker\Phi_5 &\cong \{0\}, \\
\ker\Phi_6 & \cong \integer, & \coker\Phi_6 &\cong \{0\}.
\end{align*}

\subsubsection*{Completing the proof of Theorem~\ref{3rdorientable.th}}

Putting the kernels and cokernels calculated above together using
equation~\eqref{directsum.eq} gives the table in 
Theorem~\ref{3rdorientable.th}. Taking alternating sums of Betti numbers
gives
\begin{align}
\chi(\expsig{3}) &= 3 - 4g +\binom{2g}{2}-\binom{2g}{3} \nonumber \\
                    &= \frac{-4g^3 + 12g^2 -17g + 9}{3} \label{euler.eq}
\end{align}
for $g\geq 2$, and direct substitution shows it holds for $g\geq 0$ also.
As a check we calculate the Euler characteristic using
\[
\chi(\model)=\chi(\symthree)+\chi(\symtwo)-\chi(\Sigma\times\Sigma).
\]
Macdonald gives $\chi(\sym{n}{\Sigma})=(-1)^n\binom{2g-2}{n}$, so
\begin{align*}
\chi(\model) &= -\binom{2g-2}{3}+\binom{2g-2}{2}-(2-2g)^2 \\
             &= \frac{-4g^3 + 12g^2 -17g + 9}{3},
\end{align*}
in agreement with~\eqref{euler.eq}.

\end{document}